\documentclass[letterpaper,12pt]{article}

\usepackage{dgjournal}          

\usepackage{mathptmx}

\usepackage{graphics}

\usepackage{enumerate}

\usepackage{adjustbox}

\usepackage{hyperref}

\usepackage{amsfonts}

\usepackage[authoryear,comma,longnamesfirst,sectionbib]{natbib} 

\usepackage{bm}

\usepackage{graphicx,color}

\usepackage{amsmath}

\def\ind{\begin{picture}(9,8)
         \put(0,0){\line(1,0){9}}
         \put(3,0){\line(0,1){8}}
         \put(6,0){\line(0,1){8}}
         \end{picture}
        }
\def\nind{\begin{picture}(9,8)
         \put(0,0){\line(1,0){9}}
         \put(3,0){\line(0,1){8}}
         \put(6,0){\line(0,1){8}}
         \put(1,0){{\it /}}
         \end{picture}
    }

\newtheorem{lemma}{Lemma}

\def\Var{\text{Var}}
\def\Cov{\text{Cov}}

\def\iidsim{\stackrel{\text{ind}}{\sim}}

\def\E{\mathbb{E}}

\def\Bias{\text{Bias}}
    
\begin{document}


\begin{center}
{\bfseries\Large To Adjust or Not to Adjust? Sensitivity Analysis of $M$-Bias and Butterfly-Bias}\\
\bigskip
Peng Ding and Luke Miratrix\\
Department of Statistics, Harvard University\\
Emails: \texttt{pengding@fas.harvard.edu} and \texttt{lmiratrix@stat.harvard.edu}
\end{center}

\bigskip 

\begin{center}
{\bfseries Abstract}
\end{center}

``$M$-Bias,'' as it is called in the epidemiologic literature, is the bias introduced by conditioning on a pretreatment covariate due to a particular ``$M$-Structure'' between two latent factors, an observed treatment, an outcome, and a ``collider.''
This potential source of bias, which can occur even when the treatment and the outcome are not confounded, has been a source of considerable controversy. 
We here present formulae for identifying under which circumstances biases are inflated or reduced. 
In particular, we show that the magnitude of $M$-Bias in linear structural equation models tends to be relatively small compared to confounding bias, suggesting that it is generally not a serious concern in many applied settings. 
These theoretical results are consistent with recent empirical findings from simulation studies.
We also generalize the $M$-Bias setting (1) to allow for the correlation between the latent factors to be nonzero, and (2) to allow for the collider to be a confounder between the treatment and the outcome. 
These results demonstrate that mild deviations from the $M$-Structure tend to increase confounding bias more rapidly than $M$-Bias, suggesting that choosing \emph{to} condition on any given covariate is generally the superior choice.
As an application, we re-examine a controversial example between Professors Donald Rubin and Judea Pearl.


\newpage

\section{Introduction}
The hallmark of an observational study is selection bias \citep{Heckman1979, Copas1997, Hernan2004}.
Many statisticians believe that ``there is no reason to avoid adjustment for a variable describing subjects before treatment'' in observational studies \cite[pp 76]{Rosenbaum2002}, because ``typically, the more conditional an assumption, the more generally acceptable it is'' \citep{Rubin2009}.
This advice, recently dubbed the ``pretreatment criterion'' \citep{VanderWeele2011}, is widely used in empirical studies, as more covariates generally seem to make the ignorability assumption, i.e., the assumption that conditionally on the observed pretreatment covariates,  treatment assignment is independent of the potential outcomes \citep{Rosenbaum1983}, more plausible.
And, as the validity of causal inference in observational studies relies strongly on this (untestable) assumption \citep{Rosenbaum1983}, it seems reasonable to make all efforts to render it plausible.

\begin{figure}[ht]
\centering
\begin{tabular}{ccc}
\includegraphics[width = 0.3\textwidth]{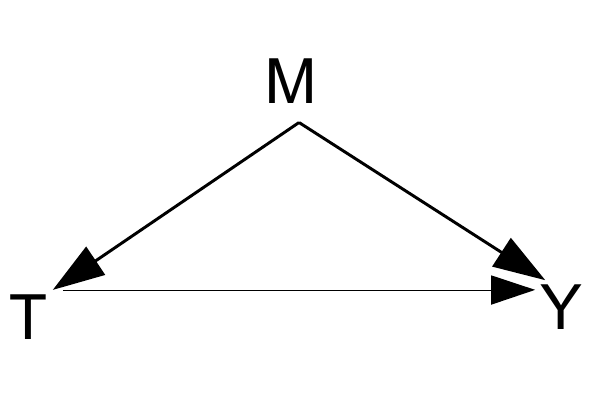}
&
\includegraphics[width = 0.3\textwidth]{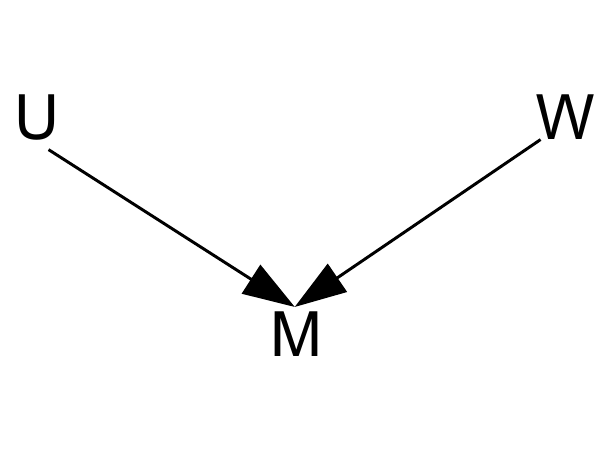}
&
\includegraphics[width = 0.3\textwidth]{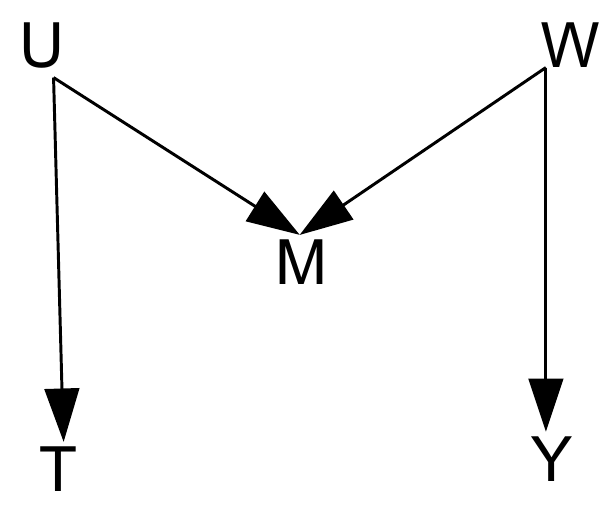}
\\
(a) A simple DAG
&
(b) V-Structure
&
(c) $M$-Structure
\\
&
 $U\ind W$ but $U\nind W|M$
&
 $T\ind Y$ but $T\nind Y|M$
\end{tabular}
\caption{Three DAGs} 
\label{fg::VandM}
\end{figure}

However, other researchers \citep{Pearl2009a, Pearl2009b, Shrier2008, Shrier2009, Sjolander2009}, mainly from the causal diagram community, do not accept this view because of the possibility of a so-called $M$-Structure, illustrated in Figure \ref{fg::VandM}(c). 
In sharp contrast to Rubin and Rosenbaum's advice, \cite{Pearl2009a} and \cite{Pearl2009b} warn practitioners that spurious bias may arise due to adjusting for a collider $M$ in an $M$-Structure, {\it even if it is a pretreatment covariate}.
This form of bias, typically called $M$-bias, a special version of so-called ``collider bias,'' has since generated considerable controversy and confusion.

We attempt to resolve some of these debates by an analysis of $M$-bias under the causal diagram or directed acyclic graph (DAG) framework.
For readers unfamiliar with the terminologies from the DAG (or Bayesian Network) literature,
more details can be found in \cite{Pearl1995} or \cite{Pearl2000}.
We here use only a small part of this larger framework.
Arguably the most important structure in the DAG, and certainly the one at root of almost all controversy, is the ``V-Structure'' illustrated in Figure \ref{fg::VandM}(b).  Here, $U$ and $W$ are marginally independent with a common outcome $M$, which shapes a ``V'' with the vertex $M$ being called a ``collider.''
From a data-generation viewpoint, one might imagine Nature generating data in two steps:  
She first picks independently two values for $U$ and $W$ from two distributions, and then she combines them (possibly along with some additional random variable) to create $M$.
Given this, conditioning on $M$ can cause a spurious correlation between $U$ and $W$, which is known as the collider bias \citep{Greenland2002}, or, in epidemiology,  Berkson's Paradox \citep{Berkson1946}.
Conceptually, this correlation happens because if one cause of an observed outcome is known to have not occurred, the other cause becomes more likely. Consider an automatic-timer sprinkler system where the sprinkler being on is independent of whether it is raining.
Here, the weather gives no information on the sprinkler.
However, given wet grass, if one observes a sunny day, one will likely conclude that the sprinklers have recently run.  Correlation has been induced.

Where things get interesting is when this collider is made into a \emph{pre-treatment} variable.
Consider Figure \ref{fg::VandM}(c), an extension of Figure \ref{fg::VandM}(b).  Here $U$ and $W$ are now also causes of the treatment $T$ and the outcome $Y$, respectively. 
Nature, as a last, third step generates $T$ as a function of $U$ and some randomness, and $Y$ as a function of $W$ and some randomness.
This structure is typically used to represent a circumstance where a researcher observes $T$, $Y$, and $M$ in nature and is attempting to derive the causal impact of $T$ on $Y$. $U$ and $W$ are unobserved, or latent.
Clearly, the causal effect of $T$ on $Y$ is zero, which is also equal to the marginal association between $T$ and $Y$.  
If a researcher regressed $Y$ on $T$, he or she would obtain a zero in expectation, which is correct for estimating the causal effect.
But perhaps there is a concern that $M$, a pretreatment covariate, may be a confounder that is masking a treatment effect.
Typically, one would then ``adjust'' for $M$ to take this possibility into account, e.g., by including $M$ in a regression or by matching units on similar values of $M$.
If we do this in this circumstance, however, then we will not find a zero causal effect, in expectation.
This is the so-called ``$M$-Bias,'' and this special structure is called the ``$M$-Structure'' in the DAG literature.

Previous qualitative analysis for binary variables shows that collider bias generally tends to be small \citep{Greenland2002},
and simulation studies \citep{Liu2012} again demonstrate that $M$-Bias is small in many realistic settings.
While mathematically describing the magnitudes of $M$-Bias in general models is intractable,
it is possible to derive exact formulae of the biases as functions of the correlation coefficients in linear structural equation models (LSEMs).
The LSEM has a long history in statistics \citep{Wright1921, Wright1934} to describe dependence among multiple random variables.
\cite{Spirtes2002} uses linear models to illustrate $M$-Bias in observational studies, and \cite{Pearl2013b} also utilize the transparency of such linear models to examine various types of causal phenomena, biases, and paradoxes. 
We here extend these works and provide exact formulae for biases, allowing for a more detailed quantitative analysis of $M$-bias. 

While $M$-Bias does exist when the true underlying data generating process (DGP) follows the exact $M$-Structure, it might be rather sensitive to various deviations from the exact $M$-Structure.
Furthermore, some might argue that an exact $M$-Structure is unlikely to hold in practice.
\cite{Gelman2011}, for example, doubts the exact independence assumption required for the $M$-Structure in the social sciences by arguing that there are ``(almost) no true zeros'' in this discipline.
Indeed, since $U$ and $W$ are often latent characteristics of the same individual, the independence assumption $U\ind W$ is a rather strong structural assumption.
Furthermore, it might be plausible that the pretreatment covariate $M$ is also a confounder between, i.e., has some causal impact on both, the treatment and outcome.
We extend our work by accounting for these departures from a pure $M$-Structure, and find that even slight departures from the $M$-Structure can 
dramatically change the forms of the biases.

This paper theoretically compares the bias from conditioning on an $M$ to not under several scenarios and finds that $M$-Bias is indeed small relative to other concerns unless there is a strong correlation structure for the variables.
We further show that these findings extend to a binary treatment regime as well.
This argument proceeds in several stages.
First, in Section \ref{sec::gaussian-m-butterfly}, we examine a pure $M$-Structure and introduce our LSEM framework.  We then discuss the cases when the latent variables $U$ and $W$ may be correlated and $M$ may also be a confounder between the treatment $T$ and the outcome $Y.$
In Section \ref{sec::binarytreatment}, we generalize the results in Section \ref{sec::gaussian-m-butterfly} to a binary treatment.
In Section \ref{sec::example-rubin-pearl}, we illustrate the theoretical findings using a controversial example between Professors Donald Rubin and Judea Pearl \citep{Rubin2007, Pearl2009b}.
Section \ref{sec::further} discusses the relevance of our findings by examining $M$-Bias in actual practice and by comparing asymptotic to finite sample properties. 
We conclude with a brief discussion and present all technical details in the Appendix.

\section{$M$-Bias and Butterfly-Bias in LSEMs}\label{sec::gaussian-m-butterfly}

We begin by examining pure $M$-Bias in a LSEM.
As our primary focus is bias, we assume data are ample and that anything estimable is estimated with nearly perfect precision.
In particular, when we say we obtain a result from a regression, we implicitly mean we obtain that result in expectation; in practice an estimator will be near the given quantities.
We do not compare relative uncertainties of different estimators given the need to estimate more or fewer parameters.
There are likely degrees-of-freedom issues that would implicitly advocate using estimators with fewer parameters, but in the circumstances considered here these concerns are likely to be minor as all the models have few parameters.

A causal DAG can be viewed as a hierarchical DGP.  In particular, any variable on the graph can be viewed as a function of its parents and some additional noise, i.e., if $R$ had parents $A, B,$ and $C$, we would have
$$
R = f( A, B, C, \epsilon_R ) \mbox{ with } \epsilon_R \ind (A,B,C).
$$
Generally noise terms such as $\epsilon_R$ are considered to be independent from each other, but they can also be given an unknown correlation structure corresponding to earlier variables not explicitly included in the diagram.  This is typically represented by drawing the dependent noise terms jointly from some multivariate distribution. 
This framework is quite general; we can represent any distribution that can be factored as a product of conditional distributions corresponding to a DAG (which is one representation of the Markov Condition, a fundamental assumption for DAGs).

LSEMs are special cases of the above with additional linearity and additivity constraints.
For simplicity, and without loss of generality, we also rescale all primary variables  $(U,W,M,T,Y)$ to have zero mean and unit variance.
For example, consider this data generating process corresponding to Figure~\ref{fg::VandM}(a):
\begin{eqnarray*}
\left\{
\begin{array}{lll}
M, \varepsilon_T, \varepsilon_Y  &\iidsim& [0, 1], \\
T &=& a M +  \sqrt{1 - a^2 } \varepsilon_T, \\
Y &=& b T + c M + \sqrt{1 - b^2 - c^2 } \varepsilon_Y,
\end{array}
\right.
\end{eqnarray*}
where we use $A\sim [0,1]$ to denote a random variable with mean zero and variance one.

In the causal DAG literature, we think about causality as reaching in and fixing a given node to a set value, but letting Nature take her course otherwise.  
For example, if we were able to set $T$ at  $t$, the above data generation process would be transformed to: 
\begin{eqnarray*}
\left\{
\begin{array}{lll}
M, \varepsilon_T, \varepsilon_Y  &\iidsim& [0,1], \\
T &=& t, \\
Y  &=& b t + c M + \sqrt{1 - b^2 - c^2 } \varepsilon_Y.
\end{array}
\right.
\end{eqnarray*}
The previous cause, $M$, of $T$ has been broken, but the impact of $T$ on $Y$ remains intact.  This changes the distribution of $Y$ but not $M$.
More importantly, this results in a distribution distinct from that of \emph{conditioning} on $T=t$.
Consider the case of positive $a,b,$ and $c$.  If we \emph{observe} a high $T$, we can infer a high $M$ (as $T$ and $M$ are correlated) and a high $Y$ due to both the $bT$ and $cM$ terms in $Y$'s equation.  However, if we \emph{set} $T$ to a high value, $M$ is unchanged.  Thus, while we will still have the large $bT$ term for $Y$, the $cM$ term will be 0 in expectation.  Thus, the expected value for $Y$ will be less.

This setting as compared to conditioning is represented with the ``do'' operator.
Given the ``do'' operator, we define a local causal effect of $T$ on $Y$ at $T=t$ as:
$$
\tau_{t} =  \frac{ \partial \E\{ Y\mid \text{do}(T) = t \} }{  \partial t } .
$$
For linear models, the local causal effect is a constant, and thus we do not need to specify $t$.
We use ``do'' here purely to indicate the different distributions. For a more technical overview, see \cite{Pearl1995} or \cite{Pearl2000}.  Our results, with more formality, can easily be expressed in this more technical notation.

\begin{figure}[ht]
\centering
\includegraphics[width = 0.4\textwidth]{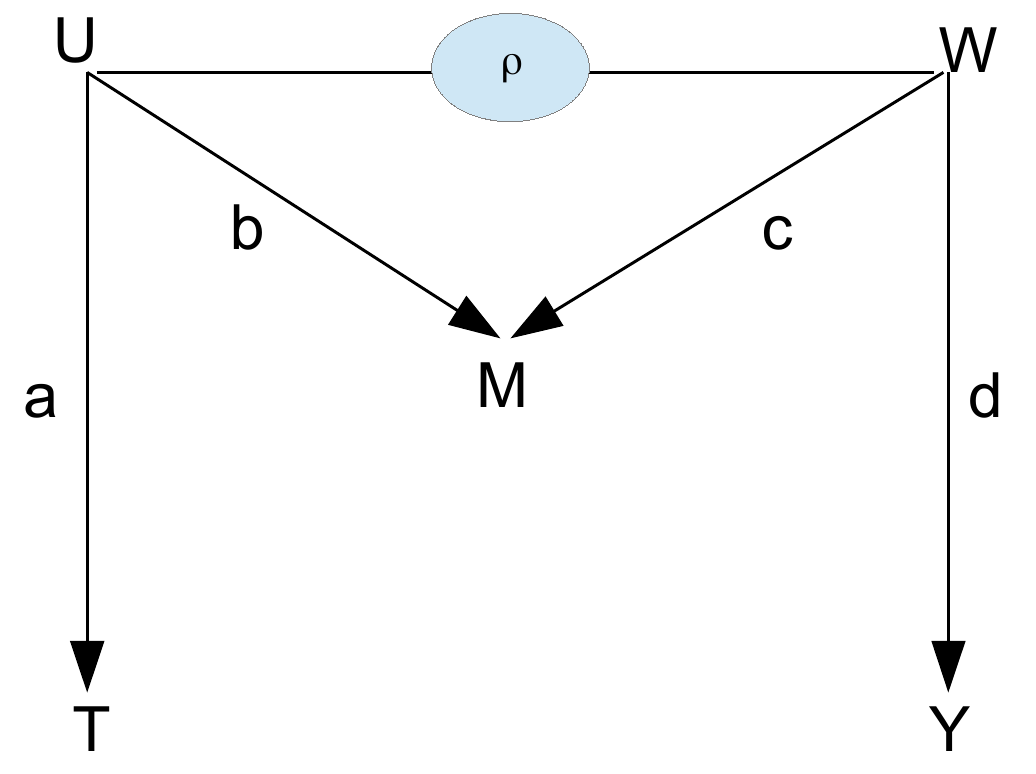}
\caption{$M$-Structure with Possibly Correlated Hidden Causes}\label{fg::mbias}
\end{figure}

If we extend the  $M$-Structure in Figure \ref{fg::VandM}(c) by allowing possible correlation between the two hidden causes $U$ and $W$, we obtain the DAG in Figure \ref{fg::mbias}. 
This in turn gives the following DGP:
\begin{eqnarray*}
\left\{
\begin{array}{lll}
\varepsilon_M, \varepsilon_T, \varepsilon_Y  &\iidsim& [0,1], \\
(U, W) &\sim& [0,0;1,1,\rho], \\
M &=& b U + cW + \sqrt{1 - b^2 - c^2} \varepsilon_M, \\
T  &=& a U + \sqrt{1- a^2} \varepsilon_T, \\
Y  &=& dW + \sqrt{1 - d^2} \varepsilon_Y,
\end{array}
\right.
\end{eqnarray*}
where we use $(A,B)\sim [0,0;1,1,\rho]$ to denote a bivariate random vector with means zero, variances one and correlation coefficient $\rho.$

Here, the true causal effect of $T$ on $Y$ is zero, namely, $\tau_t = 0$ for all $t$.
The unadjusted estimator for the causal effect obtained by regressing $Y$ onto $T$ is the same as the covariance between $T$ and $Y$:
\begin{eqnarray*}
\Bias_{unadj} = \Cov(T, Y) = \Cov(aU, dW)= ad\Cov(U, W)= ad\rho. \label{eq::bias1}
\end{eqnarray*} 
The adjusted estimator (see Lemma \ref{lemma::mbias-gaussian} in Appendix A for a proof) obtained by regressing $Y$ onto $(T, M)$ is
\begin{eqnarray*}
\Bias_{adj} =  \frac{  ad\rho(1-b^2-c^2-bc\rho) -  abcd}{  1 - (ab + ac\rho)^2    }.\label{eq::bias2}
\end{eqnarray*}
The results above and some of the results discussed later in this paper can be obtained directly from  traditional path analysis \citep{Wright1921, Wright1934}. 
However, we provide elementary proofs, which can easily be extended to binary treatment, in the Appendix.
If we allowed for a treatment effect, our results would remain essentially unchanged; the only difference would be due to restrictions on the correlation terms needed to maintain unit variance for all variables.

The above can also be expressed in the potential outcomes framework \citep{Neyman1923, Rubin1974}.  In particular, for a given unit let Nature draw $\varepsilon_M, \varepsilon_T, \varepsilon_Y, U,$ and $V$ as before.  Let $T$ be the ``natural treatment'' for that unit, i.e., what treatment it would receive sans intervention.  
Then calculate $Y(t)$ for any $t$ of interest using the ``do'' operator.  These are what we would see if we set $T=t$.  How $Y(t)$ changes for a particular unit defines that unit's collection of potential outcomes.  Then $\E\{ Y(t) \}$ for some $t$ is the expected potential outcome over the population for a particular $t$.  We can examine the derivative of this function as above to get a local treatment effect.
This connection is exact: the findings in this paper are the same as what one would find using this DGP and the potential outcomes framework.  
We here examine regression as the estimator. Note that matching would produce identical results as the amount of data grew (assuming the data generating process ensures common support, etc.).

\paragraph{Exact $M$-Bias.}  
The $M$-Bias originally considered in the literature is the special case where the correlation coefficient between $U$ and $W$ is $\rho=0$.  
In this case, the unadjusted estimator is unbiased and the absolute bias of the adjusted estimator is $|abcd|/\{ 1 - (ab)^2  \}$.
With moderate correlation coefficients $a,b,c,d$ the denominator $1-(ab)^2$ is close to one, and the bias is close to $-abcd$. 
Since $abcd$ is a product of four correlation coefficients, it can be viewed as a ``higher order bias.'' 
For example, if $a=b=c=d=0.2$, then $1-(ab)^2 = 0.9984\approx 1$, and the bias of the adjusted estimator is $ -abcd/\{ 1 - (ab)^2  \} = -0.0016\approx 0$;
if $a=b=c=d=0.3$, then $1-(ab)^2 = 0.9919\approx 1$, and the bias of the adjusted estimator is $- abcd/\{ 1 - (ab)^2  \} = - 0.0082 \approx 0.$  
Even moderate correlation results in little bias.

\begin{figure}[ht]
\centering
\includegraphics[width = \textwidth]{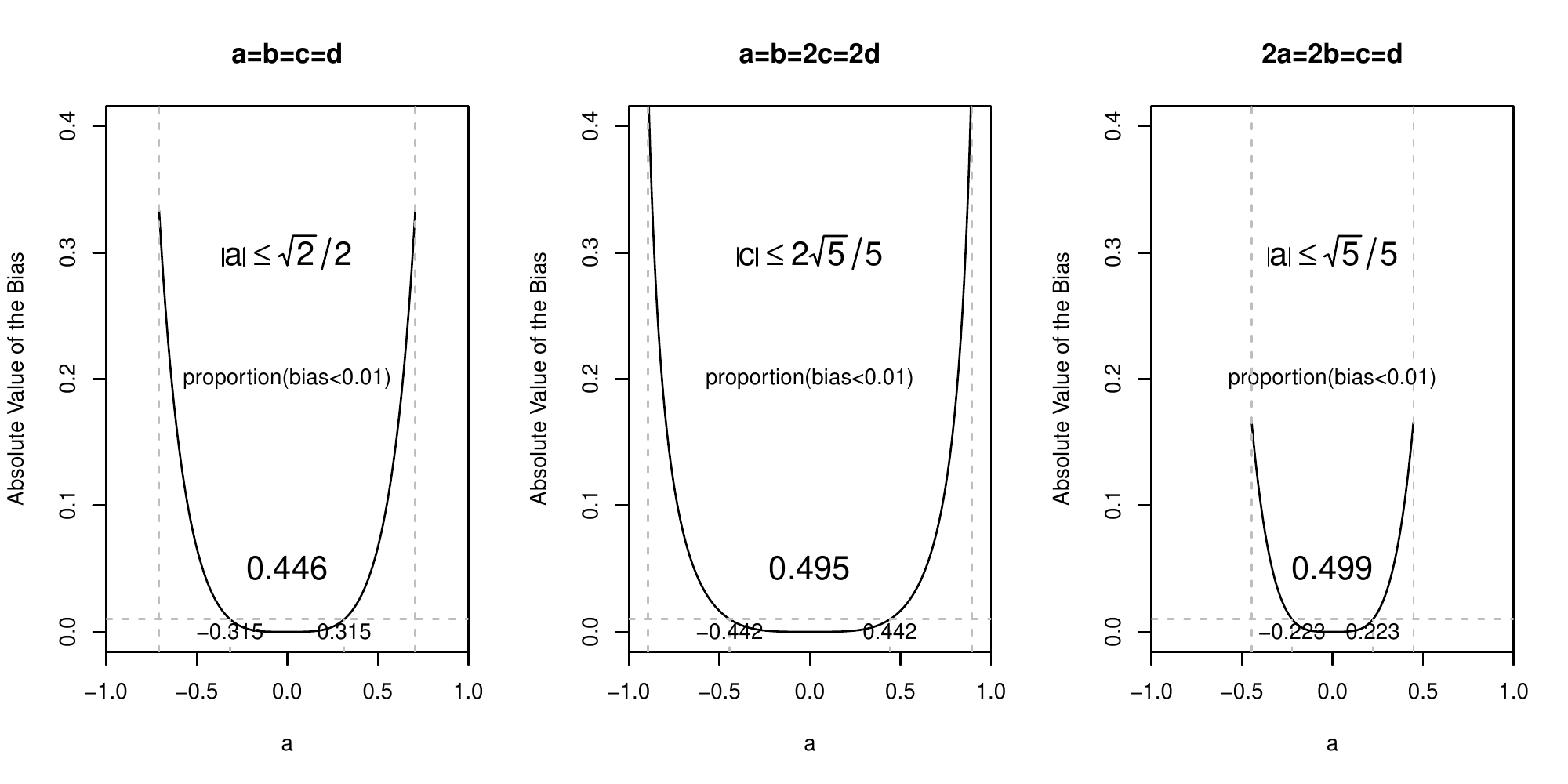}

\caption{$M$-Bias with Independent $U$ and $W$.  The three subfigures correspond to the cases when $(U,W)$ are equally/more/less predictive to the treatment than to the outcome.
In each subfigure, we show the proportions of the areas where the adjusted estimator has a bias smaller than $0.01$.}\label{fg::mbias-gaussian-ind}
\end{figure}

In Figure~\ref{fg::mbias-gaussian-ind}, we plot the bias of the adjusted estimator as a function of the correlation coefficients, and let these coefficients change to see how the bias changes.
In the first subfigure, we assume all the correlation coefficients have the same magnitude ($a=b=c=d$), and we plot the absolute bias of the adjusted estimator versus $a$. 
The constraints on variance and correlation only allow for some combinations of values for $a, b, c$ and $d$ which limits the domain of the figures.
In this case, for example, $|a|\leq \sqrt{2}/2$ due to the requirement that $b^2+c^2 = 2a^2 \leq 1$.  Other figures have limited domains due to similar constraints.
In the second subfigure of Figure~\ref{fg::mbias-gaussian-ind}, we assume that $M$ is more predictive to the treatment $T$ than to the outcome $Y$, with $a=b=2c=2d$.
In the third subfigure of Figure~\ref{fg::mbias-gaussian-ind}, we assume that $M$ is more predictive to the outcome $Y$, with $2a=2b=c=d$.
The biases are generally very small within wide ranges of the feasible regions of the correlation coefficients. 
However, the biases do blow up when the correlation coefficients are extremely large. Near the boundary of the feasible regions in Figure \ref{fg::mbias-gaussian-ind}, the $M$-Structure is approximately deterministic, which is rare in social sciences. \cite{Pearl2009a} does not exclude the worst cases, and thus he considers $M$-Bias as a severe problem.

\begin{figure}[ht]
\begin{tabular}{p{0.5\columnwidth} p{0.5\columnwidth} }
\includegraphics[width = 0.5\textwidth]{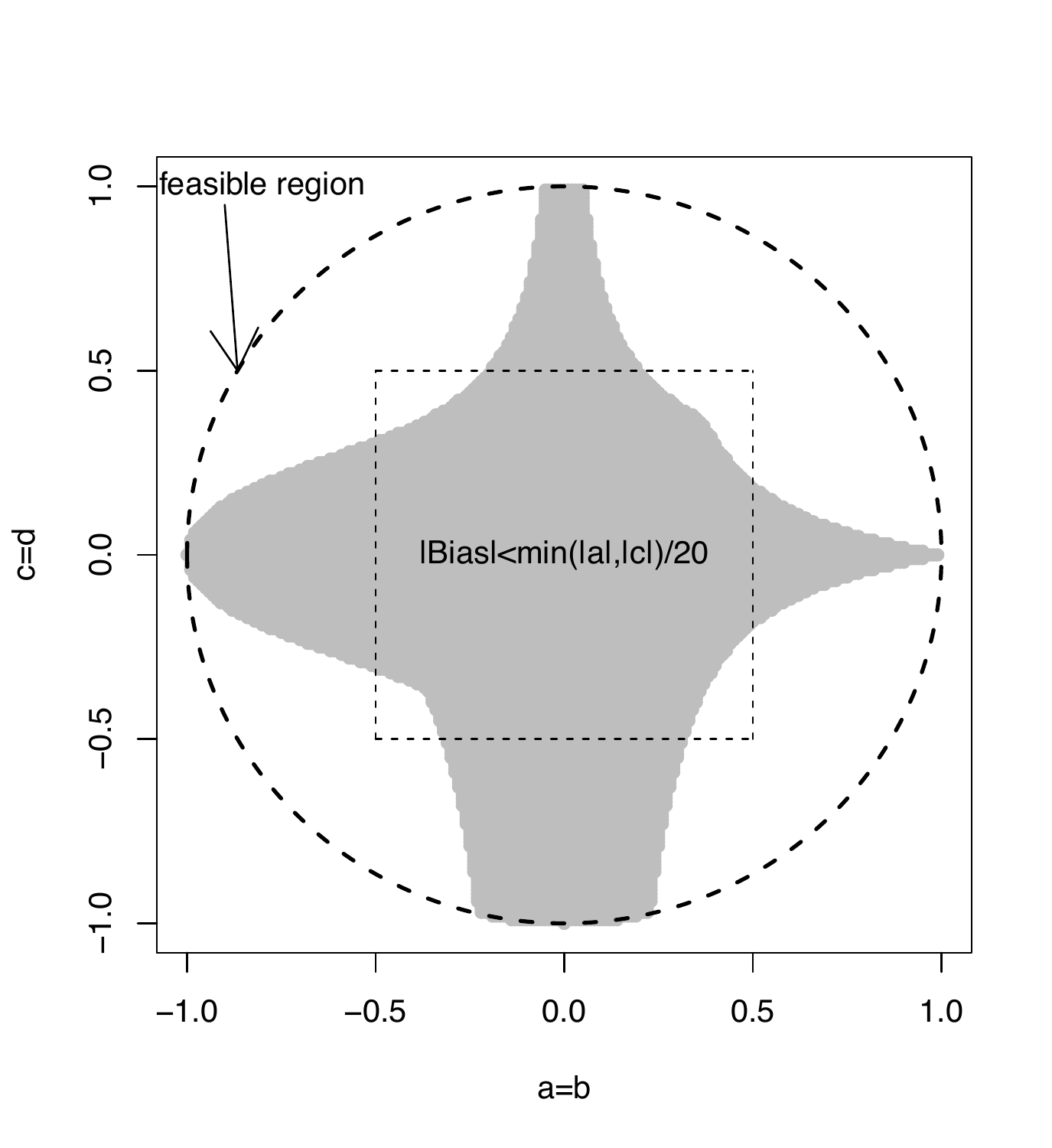}
&
\includegraphics[width = 0.5\textwidth]{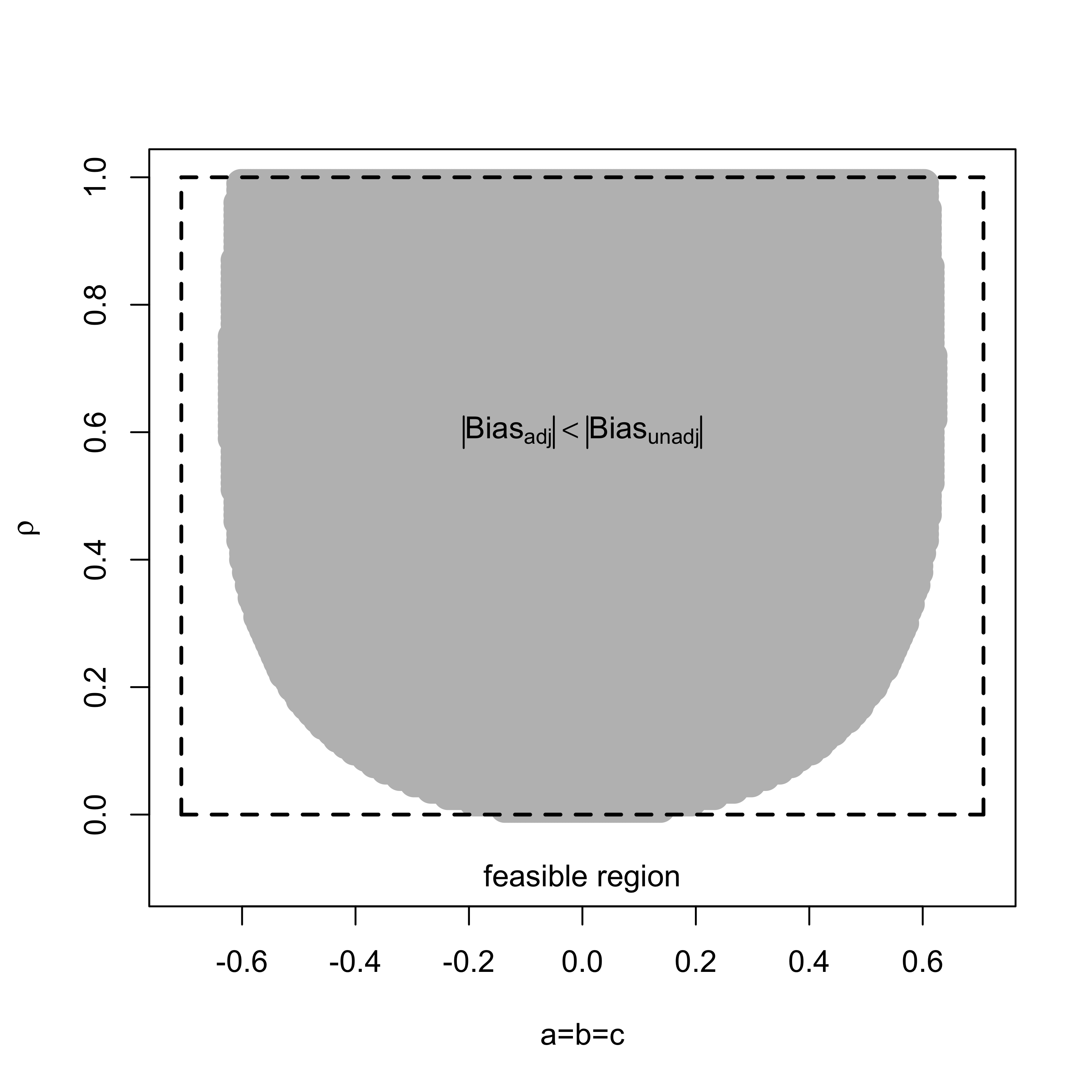}
\\
(a) Pure $M$-Bias.  Within the grey region, the absolute bias of the adjusted estimator is less than $1/20$ of the minimum of $|a|(=|b|)$ and $|c|(=|d|)$.
&
(b) $M$-Bias with Correlated $U$ and $W$.  Within the grey region, the adjusted estimator is superior. 
\end{tabular}
\caption{$M$-Bias under Different Scenarios}\label{fg::mbias-grey-plots}
\end{figure}

In Figure~\ref{fg::mbias-grey-plots}(a), we assume $a=b$ and $c=d$ and examine a broader range of relationships.  Here, the grey area satisfies $|\text{Bias}_{adj}| < \min(|a|,|c|)/20$.
For example, when the absolute values of the correlation coefficients are smaller than $0.5$ (the square with dashed boundary in Figure~\ref{fg::mbias-grey-plots}(a)), the corresponding area is almost grey, implying small bias.

Due to the four dimensional sensitivity parameters $(a,b,c,d)$, a full exploration and graphical illustration over all possible values of the sensitivity parameters is formidably hard. In the absence of prior knowledge about the DAG, our sensitivity analysis here is based on some simplifications (e.g., $a=b$ and $c=d$), which may reflect some real situations. Using the bias formulae in this paper, we can easily conduct sensitivity analysis for other parameter combinations, depending on our practical problem and background knowledge about the DAG.

As a side note, \cite{Pearl2013b} noticed a surprising fact: the stronger the correlation between $T$ and $M$, the larger the absolute bias of the adjusted estimator, since the absolute bias is monotone increasing in $|ab|.$
From the second and the third subfigure of Figure \ref{fg::mbias-gaussian-ind}, we see that when $M$ is more predictive of the treatment, the biases of the adjusted estimator indeed tends to be larger.

\paragraph{Correlated Latent Variables.} 
When the latent variables $U$ and $W$ are correlated with $\rho\neq 0$, 
\emph{both} the unadjusted and adjusted estimators may be biased.
The question then becomes: which is worse?
The ratio of the absolute biases is
$$
\left| \frac{ \text{Bias}_{adj} }{ \text{Bias}_{unadj}} \right| = \Big|     \frac{  \rho(1-b^2 - c^2-bc\rho) -   bc   }{  \rho  \{   1 - (ab+ac\rho)^2\}  }  \Big|,
$$
which does not depend on $d$ (the relationship between $W$ and $Y$).
For example, if the correlation coefficients $a,b,c,\rho$ all equal $0.2$, the ratio above is $0.714$; in this case the adjusted estimator is superior to the unadjusted one by a factor of $1.4$.
Figure~\ref{fg::mbias-grey-plots}(b) compares this ratio to 1 for all combinations of $\rho$ and $a(=b=c)$.  Generally, the adjusted estimator has smaller bias except when $a, b,$ and $c$ are quite large.

\begin{figure}[ht]
\centering
\includegraphics[width = \textwidth]{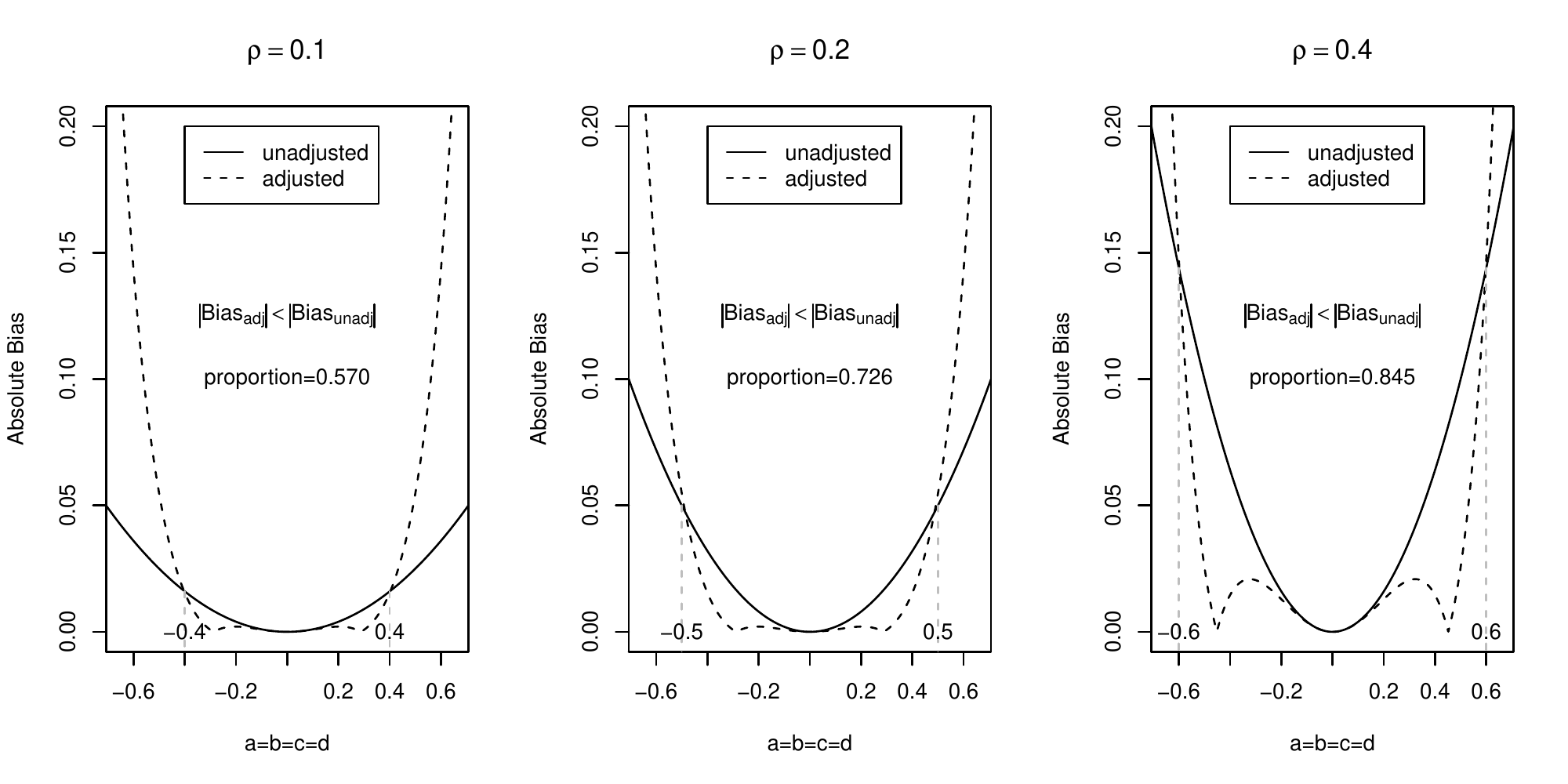}
\caption{$M$-Bias with correlated $U$ and $W$ when $\rho \neq  0$ with $a=b=c=d$. In each subfigure, we show the proportion of the areas where the adjusted estimator has a smaller bias than the unadjusted estimator.}
\label{fg::mbias-gaussian}
\end{figure}

In Figure \ref{fg::mbias-gaussian}, we again assume $a=b=c=d$ and investigate the absolute biases
as functions of $a$ 
for fixed $\rho$ at $0.1,0.2, $ and $0.4.$
When the correlation coefficients $a(=b=c=d)$ are not dramatically larger than $\rho$, the adjusted estimator has smaller bias than the unadjusted one.

\paragraph{The Disjunctive Cause Criterion.}
In order to remove biases in observational studies, \cite{VanderWeele2011} propose a new ``disjunctive cause criterion'' for selecting confounders, which requires controlling for all the covariates that are either causes of the treatment, causes of the outcome, or causes of both. According to the ``disjunctive cause criterion,'' when $\rho\neq 0$, we should control for $(U,W)$ if possible. Unfortunately, neither of $(U,W)$ is observable. However, controlling the ``proxy variable'' $M$ for $(U,W)$ may reduce bias when $\rho$ is relatively large. 
In the special case with $b=0$, the ratio of the absolute biases is
$$
\left| \frac{ \text{Bias}_{adj} }{ \text{Bias}_{unadj}} \right|  = \frac{1-c^2}{1-(ac\rho)^2} \leq 1;
$$
in another special case with $c=0$, the ratio of the absolute biases is
\begin{eqnarray}
\left| \frac{ \text{Bias}_{adj} }{ \text{Bias}_{unadj}} \right|  = \frac{1-b^2}{1-(ab)^2} \leq 1.
\label{eq::mbias-c0}
\end{eqnarray}
Therefore, if either $U$ or $W$ is not causative to $M$, the adjusted estimator is always better than the unadjusted one.

\begin{figure}[ht]
\centering
\includegraphics[width = 0.4\textwidth]{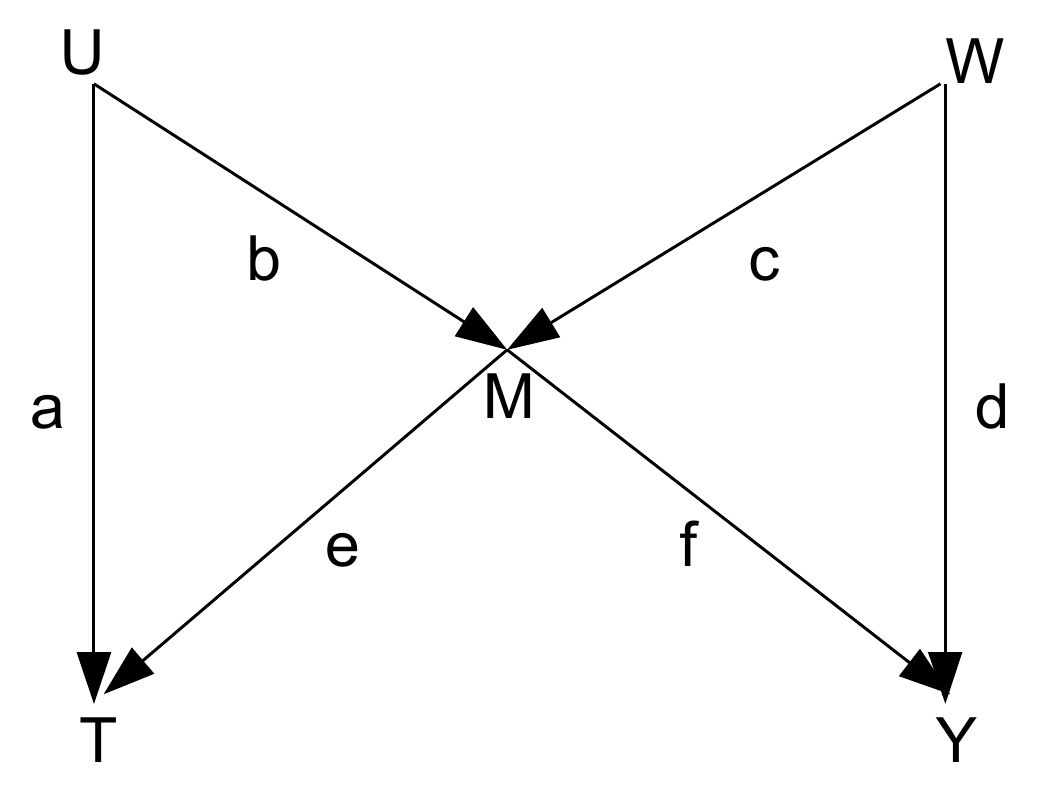}
\caption{Butterfly-Structure}\label{fg::mbias_confounding}
\end{figure}

\paragraph{Butterfly-Bias: $M$-Bias with Confounding Bias.}
Models, especially in the social sciences, are approximations. They rarely hold exactly.
In particular, for any covariate $M$ of interest, there is likely to be some concern that $M$ is indeed a confounder, even if it is also a possible source of $M$-Bias.
If we let $M$ both be a confounder as well as the middle of an $M$-Structure we obtain a ``Butterfly-Structure'' \citep{Pearl2013b} as shown in Figure \ref{fg::mbias_confounding}.
In this circumstance, conditioning will help with confounding bias, but hurt with $M$-Bias. 
Ignoring $M$ will not resolve any confounding, but will avoid $M$-Bias. 
The question then becomes that of determining which is the lesser of the two evils.

We can examine this trade-off for a LSEM corresponding to Figure \ref{fg::mbias_confounding}.  The DGP is given by the following equations:
\begin{eqnarray*}
\left\{
\begin{array}{lll}
U, W, \varepsilon_M , \varepsilon_T, \varepsilon_Y &\iidsim& [0,1],\\
M &=& b U + cW + \sqrt{1 - b^2 - c^2} \varepsilon_M,\\
T  &=& a U + eM  + \sqrt{1- a^2 - e^2} \varepsilon_T,\\
Y  &=& dW + f M + \sqrt{1 - d^2 - f^2} \varepsilon_Y.
\end{array}
\right.
\end{eqnarray*}

Again, the true causal effect of $T$ on $Y$ is zero.
The unadjusted estimator obtained by regressing $Y$ onto $T$ is the covariance between $T$ and $Y$:
$$
\Bias_{unadj} = \Cov(T, Y) = abf + cde + ef.
$$
It is not, in general, zero, implying bias.
The adjusted estimator (see Lemma \ref{lemma::butterfly-gaussian} in Appendix for a proof) obtained by regressing $Y$ onto $(T,M)$ has bias
\begin{eqnarray*}
\Bias_{adj} = -  \frac{abcd}{  1 - (ab + e)^2    }.\label{eq::bias3}
\end{eqnarray*}

If the values of $e$ and $f$ are relatively high (i.e., $M$ has a strong effect on both $T$ and $Y$), the confounding bias is large and the unadjusted estimator will be severely biased.
For example, if $a,b, c,d , e,$ and $f$ all equal $0.2$, the bias of the unadjusted estimator is $0.056$, but the bias of the adjusted estimator is only $-0.0017$, an order of magnitude smaller.
Generally, the largest term for the unadjusted bias is the second-order term of $ef$, while the adjusted bias only has, ignoring the denominator, a fourth-order term of $abcd$.  This suggests adjustment is generally preferable and that $M$-bias is in some respect a ``higher order bias.''

\begin{figure}[ht]
\begin{tabular}{p{0.5\columnwidth} p{0.5\columnwidth} }
\includegraphics[width = 0.48\textwidth]{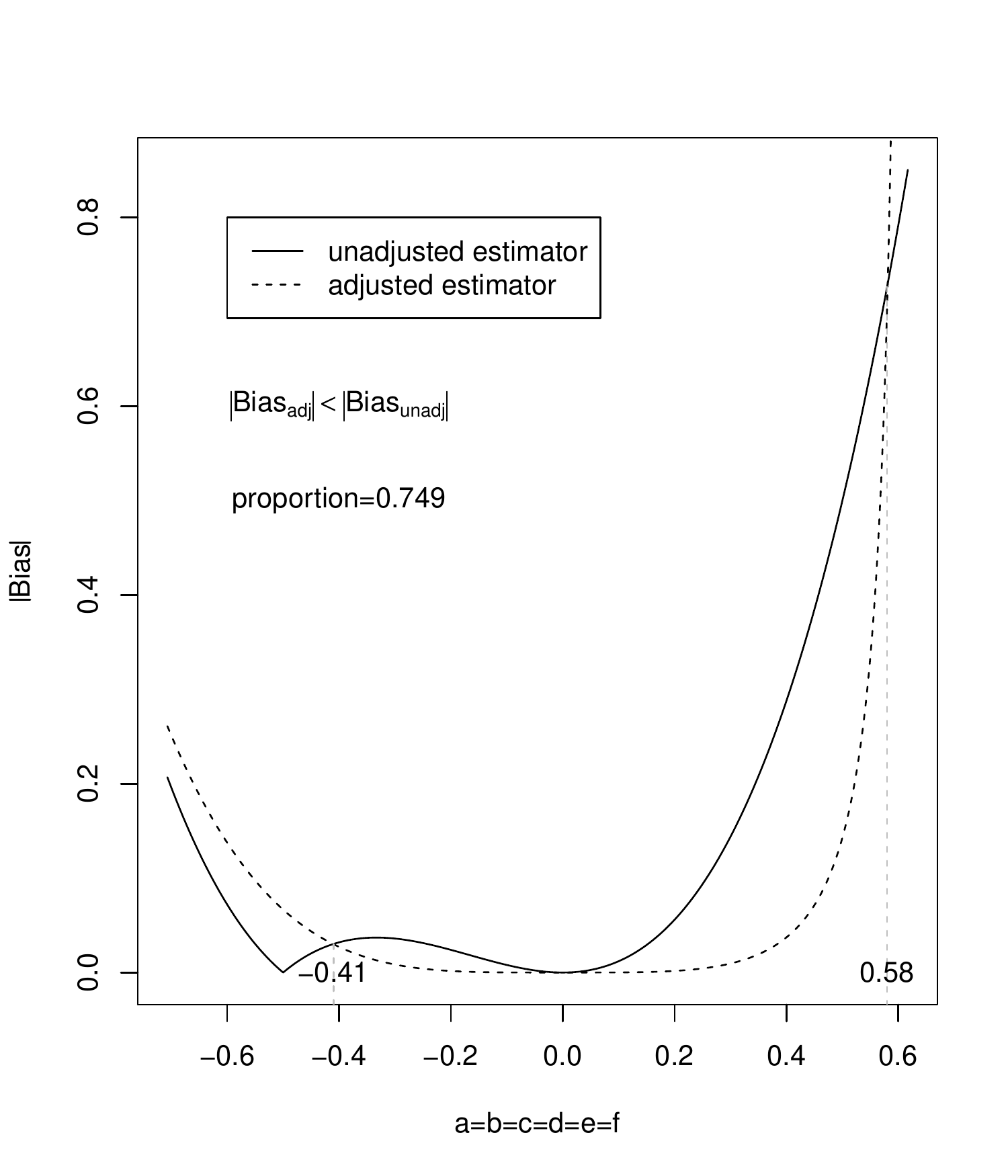}
&
\includegraphics[width = 0.48\textwidth]{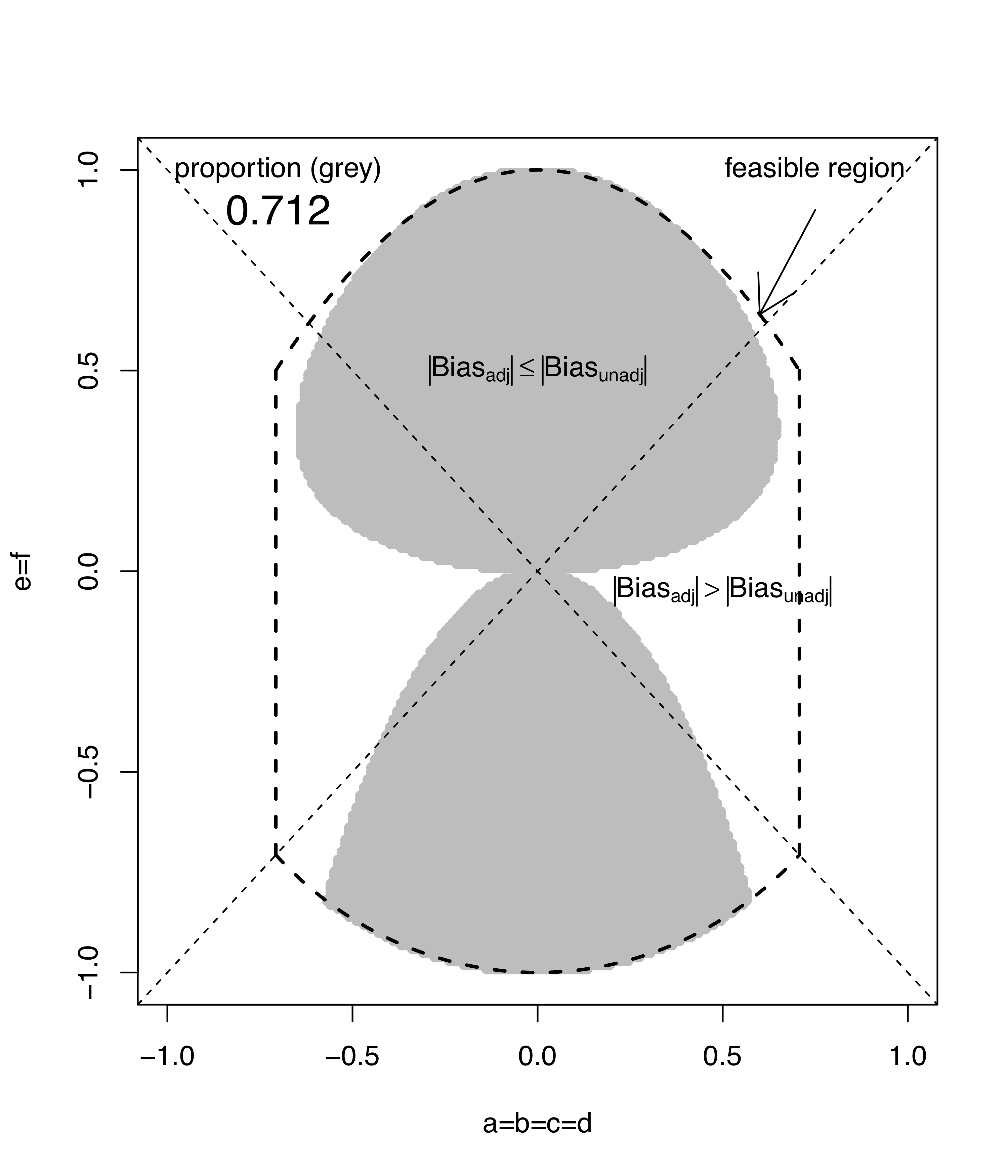}
\\
(a) Absolute biases of both estimators with $a=b=c=d=e=f$. 
&
(b) Comparison of the absolute biases with $a=b=c=d$ and $e=f$. 
Within $71.2\%$ (in grey) of the feasible region, the adjusted estimator has smaller bias than the unadjusted one.
\end{tabular}
\caption{Butterfly-Bias}\label{fg::butterfly-bias}
\end{figure}

Detailed comparison of the ratio of the biases is difficult, since we can vary six parameters $(a,b,c,d,e,f)$.
In Figure \ref{fg::butterfly-bias}(a), we assume all the correlation coefficients have the same magnitude, and plot bias for both estimators as a function of the correlation coefficient within the feasible region, defined by the restrictions
$-\sqrt{2}/2  <  a  < (-1+\sqrt{5})/2$, due to the restrictions 
\begin{eqnarray}\label{eq::restriction-butterfly}
b^2+c^2 <  1,~  a^2+e^2 <  1, ~ d^2+f^2 <  1, \text{ and } | a^2+e | < 1.
\end{eqnarray}
Within $74.9\%$ of the feasible region, the adjusted estimator has smaller bias than the unadjusted one. 
The unadjusted estimator only has smaller bias than the adjusted estimator when the correlation coefficients are extremely large.
In Figure \ref{fg::butterfly-bias}(b), we assume $a=b=c=d$ and $e=f$, and compare $|\text{Bias}_{adj} |$ and $| \text{Bias}_{unadj} |$ within the feasible region of $(a,e)$ defined by (\ref{eq::restriction-butterfly}).
We can see that the adjusted estimator is superior to the unadjusted one for $71\%$ (colored in grey in Figure \ref{fg::butterfly-bias}(b)) of the feasible region.
In the area satisfying $|e|>|a|$ in Figure \ref{fg::butterfly-bias}(b), where the connection between $M$ to $T$ and $Y$ is stronger than the other connections, the area is almost entirely grey suggesting that the adjusted estimator is preferable.  
This is sensible because here the confounding bias has larger magnitude than the $M$-Bias.
In the area satisfying $|a|<|e|$, where $M$-bias is stronger than confounding bias, the unadjusted estimator is superior for some values, but still tends to be inferior when the correlations are roughly the same size.

\section{Extensions to a Binary Treatment}\label{sec::binarytreatment}

\begin{figure}[ht]
\begin{tabular}{p{0.5\columnwidth} p{0.5\columnwidth} }
\includegraphics[width = 0.45\textwidth]{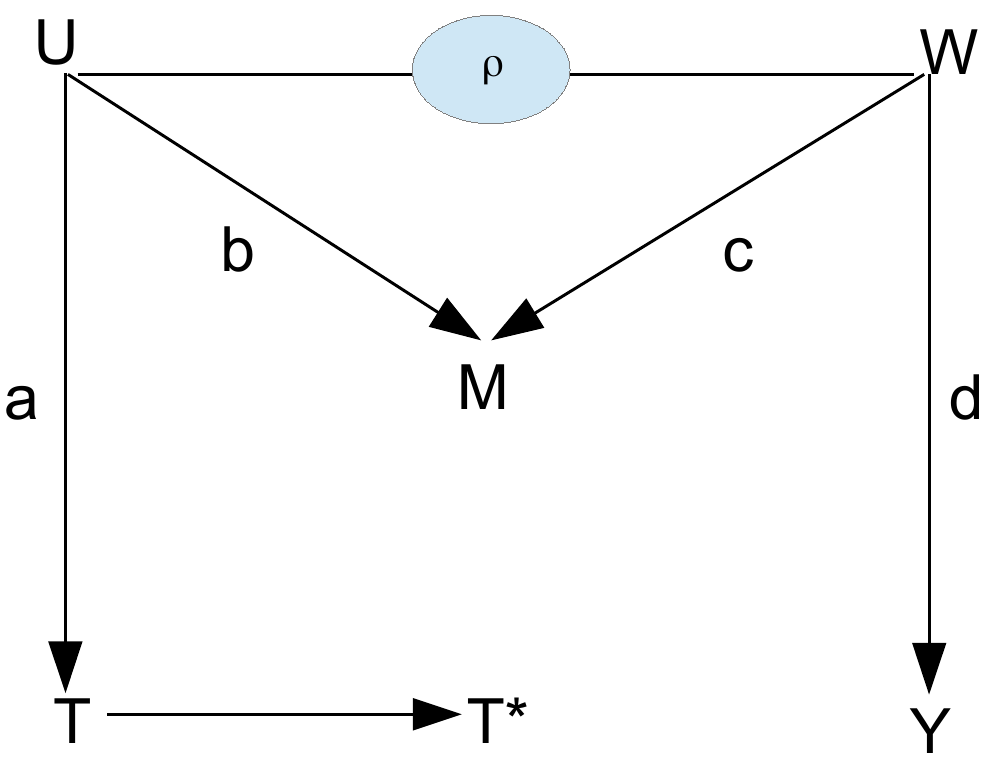}
&
\includegraphics[width = 0.45\textwidth]{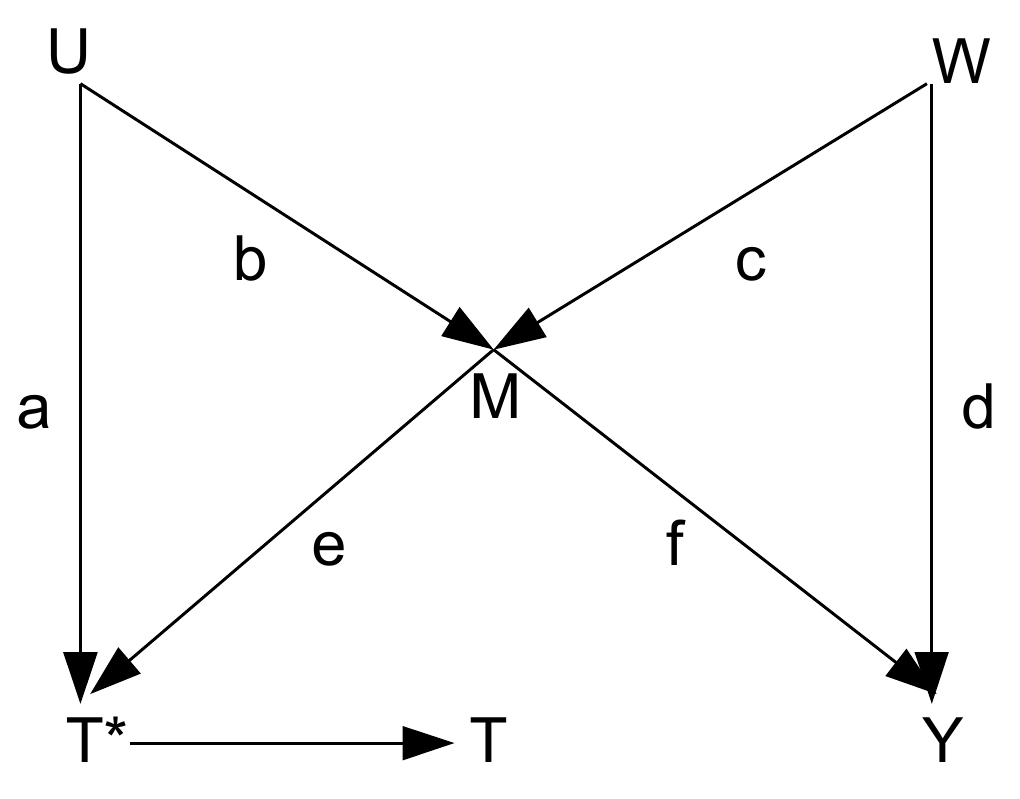}
\\
(a) Correlated Hidden Causes 
&
(b) Butterfly-Structure
\end{tabular}
\caption{DAGs for Binary Treatment}\label{fg::binary}
\end{figure}

%
%
%
%
%

One might worry that the conclusions in the previous section are not applicable for a binary treatment.
It turns out, however, that they are.
In this section, we  extend the results in Section \ref{sec::gaussian-m-butterfly} to binary treatments by representing the treatment through a latent Gaussian variable as shown in Figure \ref{fg::binary}.

\paragraph{Correlated Latent Variables.}
We extend Figure \ref{fg::mbias} to Figure \ref{fg::binary}(a).   Here, $T^*$ is the old $T$.  The generating equations for $T$ and $T^*$ become
$$
T = I(T^*  \geq  \alpha), \text{ and } T^* =  a U +\sqrt{1-a^2} \varepsilon_T.
$$
Other variables and noise terms remain the same. 
Although it might be relaxed, we make reference to the Normally assumption of the error terms for mathematical simplicity. 
The intercept $\alpha$ determines the proportion of the individuals receiving the treatment: $\Phi(-\alpha) = P(T=1)$, where $\Phi(\cdot)$ is the cumulative distribution function of a standard Normal distribution.
When $\alpha = 0$, the number of individuals exposed to the treatment and control are balanced; when $\alpha <  0$, more individuals are exposed to the treatment; when $\alpha > 0 $, the reverse.

The true causal effect of $T$ on $Y$ is again zero.
Let $\phi(\cdot) = \Phi'(\cdot)$ and $\eta(\alpha) \equiv  \phi (\alpha) /\{ \Phi(\alpha) \Phi(-\alpha) \}$.
Then Lemma \ref{lemma::3} in Appendix
shows that the unadjusted estimator has bias
\begin{eqnarray*}\label{eq::unad_bias_01}
\Bias_{unadj} = ad\rho \eta(\alpha),
\end{eqnarray*}
and the adjusted estimator has bias
$$
\Bias_{adj} = \frac{ ad\eta(\alpha) \{  \rho (1-b^2-c^2-bc\rho) - bc \}   } { \rho \{     1 - (ab + ac\rho)^2 \phi(\alpha)  \eta(\alpha)   \}   } .
$$

\begin{figure}[ht]
\centering
\includegraphics[width = \textwidth]{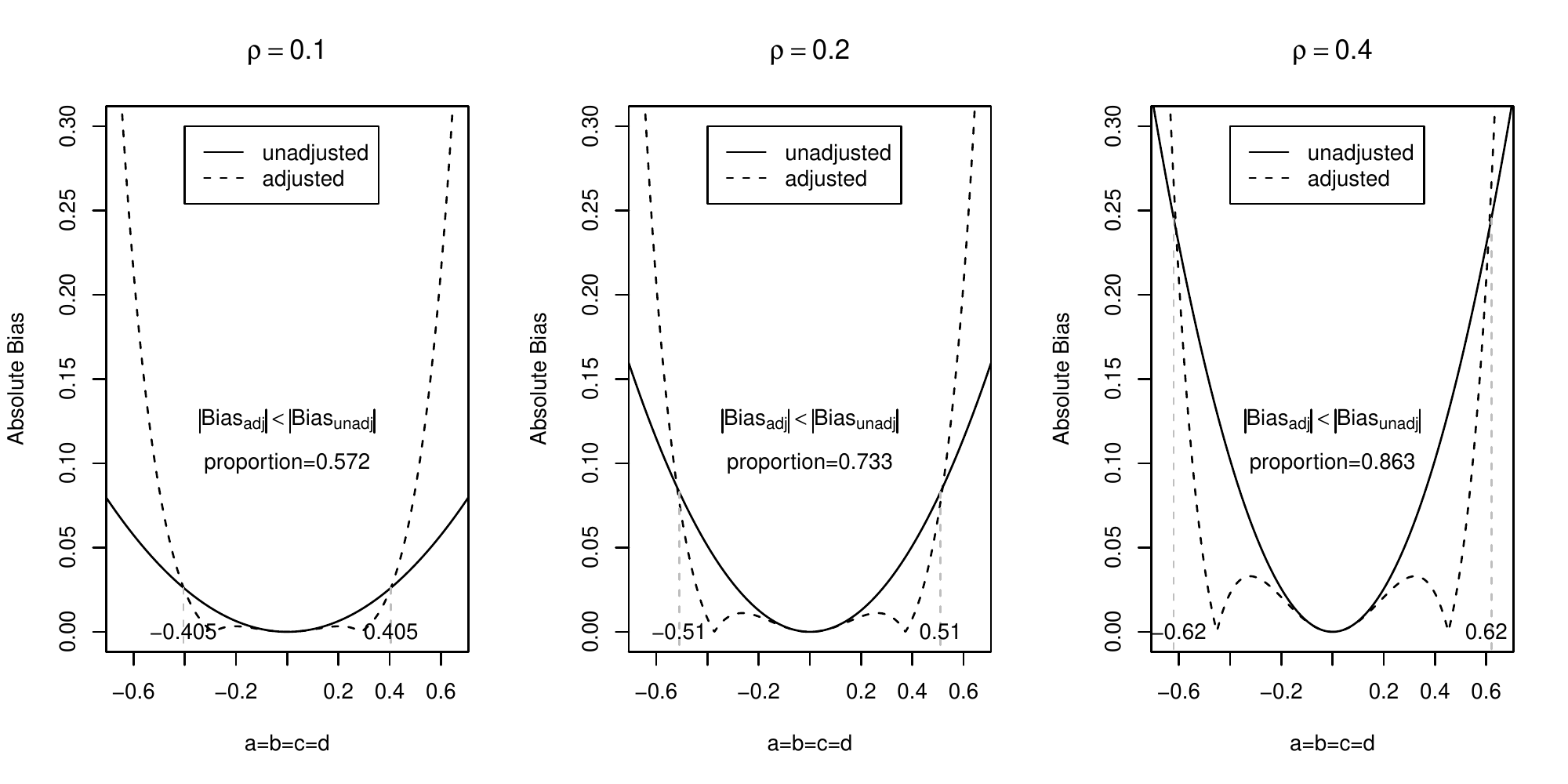}
\caption{$M$-Bias with correlated $(U,W)$ and binary treatment. Compare to Figure~\ref{fg::mbias-gaussian}}\label{fg::mbias-binary}
\end{figure}

When $\rho = 0$, the unadjusted estimator is unbiased, but the adjusted estimator has bias
$$
-\frac{abcd \eta(\alpha) }{ 1 - (ab)^2\phi(\alpha)\eta(\alpha)}.
$$
When $\rho \neq 0$, the ratio of the absolute biases is
$$
\left| \frac{ \text{Bias}_{adj} }{ \text{Bias}_{unadj}} \right|  =      \Big|    \frac{  \rho (1-b^2-c^2-bc\rho) - bc   } { \rho \{     1 - (ab + ac\rho)^2 \phi(\alpha)  \eta(\alpha)   \}   }        \Big|.
$$

The patterns for a binary treatment do not differ much from a continuous treatment.
As before, if the correlation coefficient is moderately small, the $M$-Bias also tends to be small.
As shown in Figure~\ref{fg::mbias-binary} (analogous to Figure~\ref{fg::mbias-gaussian}), when $|\rho|$ is comparable to $|a|(=|b|=|c|=|d|)$, the adjusted estimator is less biased than the unadjusted estimator. Only when $|a|$ is much larger than $|\rho|$ is the unadjusted estimator superior.

\paragraph{Butterfly-Bias with a Binary Treatment.}
We can extend the LSEM Butterfly-Bias setup to binary treatment just as we extended the $M$-Bias setup.  
Compare Figure \ref{fg::binary}(b) to Figure \ref{fg::mbias_confounding}.
$T$ becomes $T^*$ and $T$ is built from $T^*$ as above.
The structural equations for $T$ and $T^*$ for butterfly bias in the binary case are then
$$
T = I(T^*  \geq  \alpha ), \text{ and } T^* =   a U + eM +\sqrt{1-a^2 - e^2}  \varepsilon_T.
$$
The other equations and variables are the same as before.


Although the true causal effect of $T$ on $Y$ is zero,
Lemma \ref{lemma::4} in Appendix shows that
the unadjusted estimator has bias
\begin{eqnarray} \label{eq::butterfly-binaryT-unadjust}
\Bias_{unadj} = (cde + abf + ef) \eta(\alpha),
\end{eqnarray}
and the adjusted estimator has bias
\begin{eqnarray}
\Bias_{adj} = - \frac{  abcd \eta(\alpha)  }{ 1 -  (ab + e)^2\phi(\alpha)\eta(\alpha)}.
\end{eqnarray}
Therefore, the ratio of the absolute biases is
$$
\left| \frac{ \text{Bias}_{adj} }{ \text{Bias}_{unadj}} \right| = 
\Big|     
\frac{  abcd \eta(\alpha) }
{   (cde + abf +ef) \{1 - (ab + e)^2 \phi (\alpha)\eta(\alpha)  \}   } \Big| .
$$

\begin{figure}[ht]
\begin{tabular}{p{0.5\columnwidth} p{0.5\columnwidth}}
\includegraphics[width = 0.48\textwidth]{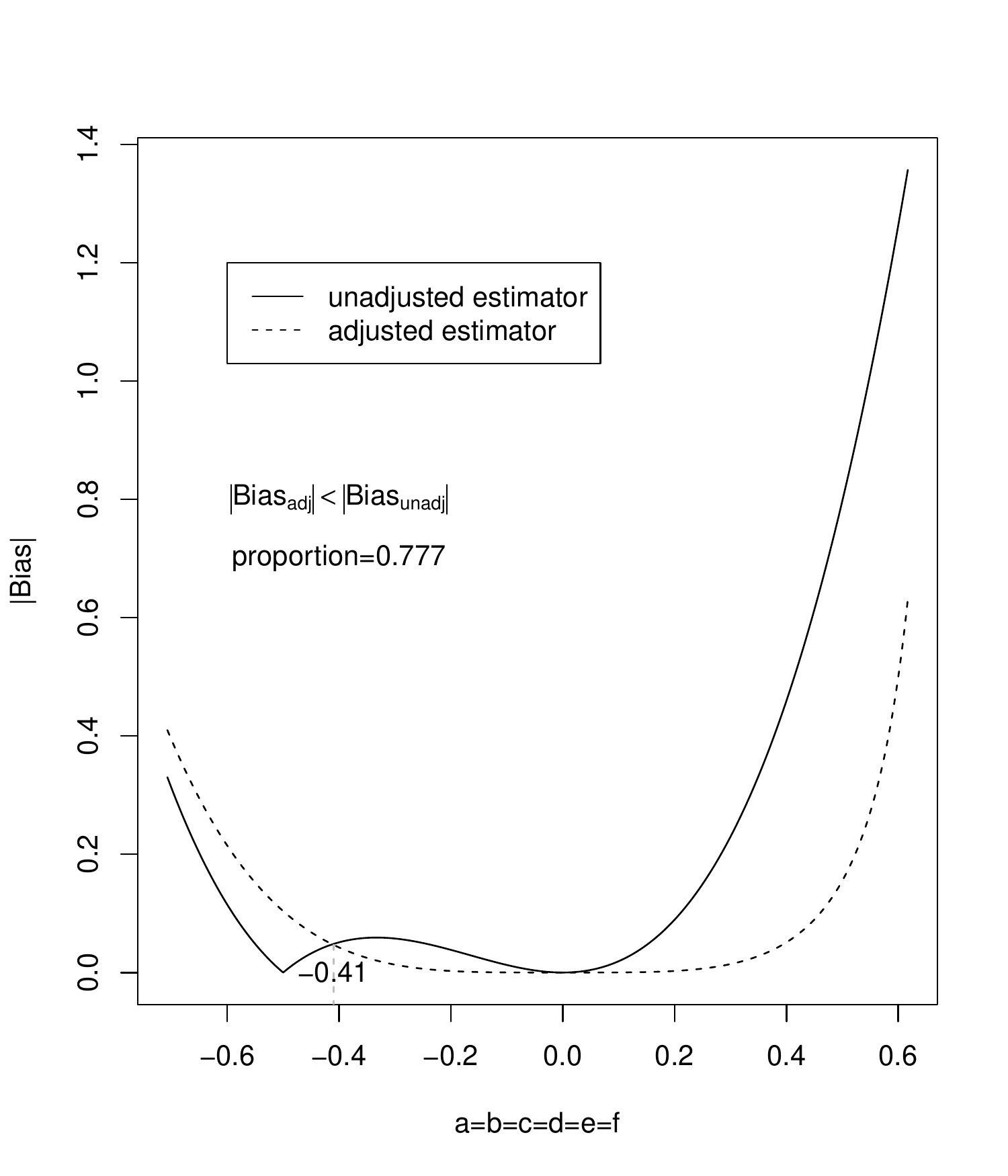}
&
\includegraphics[width = 0.48\textwidth]{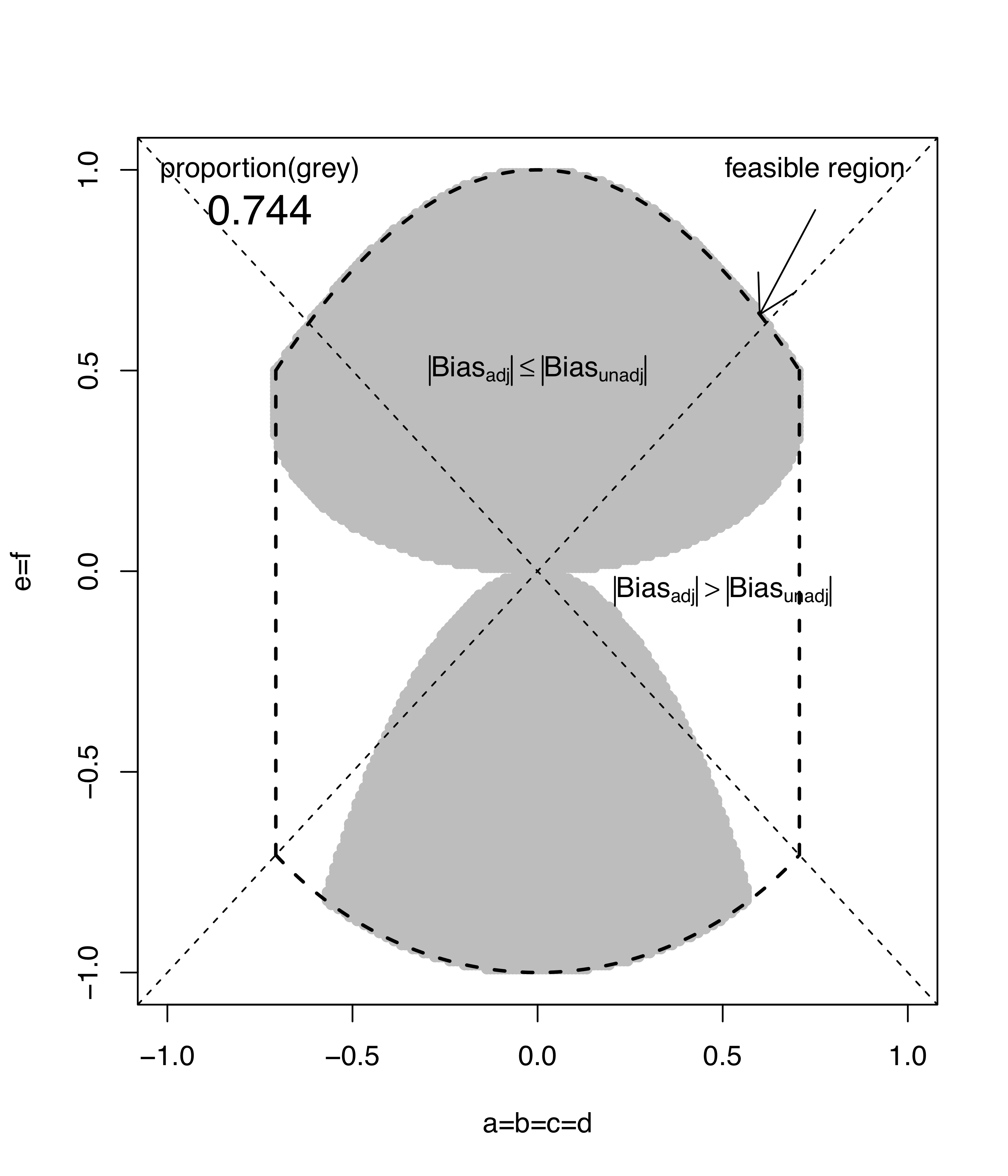}
\\
(a) Absolute biases with $a=b=c=d=e=f$. 
&
(b) Comparison of the absolute bias with $a=b=c=d$ and $e=f$. Within $74.4\%$ (in grey) of the feasible region, the adjusted estimator is better than the unadjusted estimator.
\end{tabular}
\caption{Butterfly-Bias with a Binary Treatment}\label{fg::binary-butterfly-bias}
\end{figure}

Complete investigation of the ratio of the biases is intractable with seven varying parameters $(a,b,c,d,e,f,\alpha)$.
However, in the very common case with $\alpha=0$, which gives equal-sized treatment and control groups, we again find trends similar to the continuous treatment case.
See Figure~\ref{fg::binary-butterfly-bias}.
As before, only in the cases with very small $e (= f)$ but large $a(=b=c=d)$, does the unadjusted estimator tend to be superior.
Within a reasonable region of $\alpha$, these patterns are quite similar.

\section{Illustration: The Rubin--Pearl Controversy}\label{sec::example-rubin-pearl}

\begin{figure}[ht]
\centering
\footnotesize
\begin{tabular}{p{0.5\columnwidth} p{0.5\columnwidth} }
\adjustbox{valign=m}{\includegraphics[width = 0.48\textwidth]{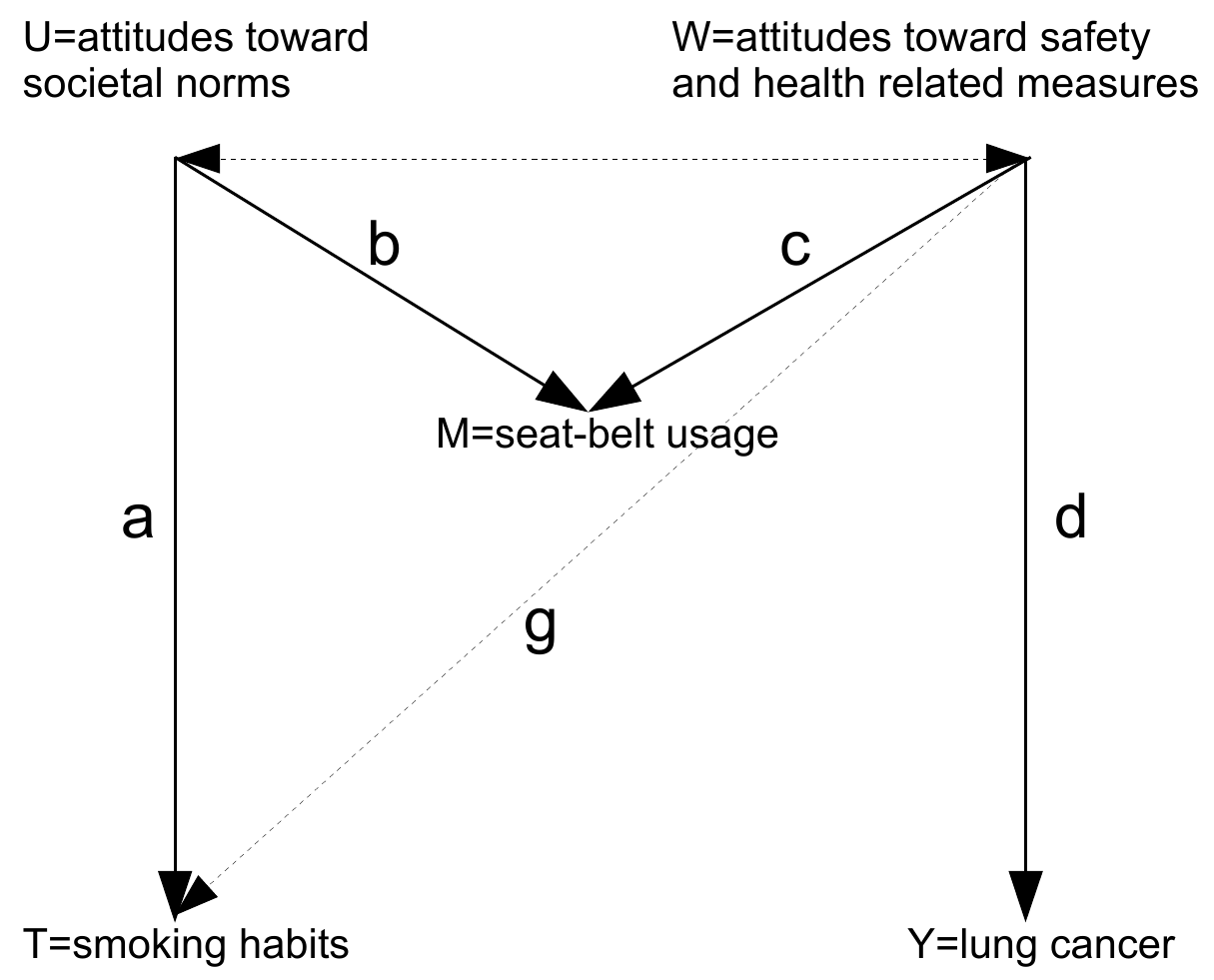}}
&
\adjustbox{valign=m}{\includegraphics[width = 0.48\textwidth]{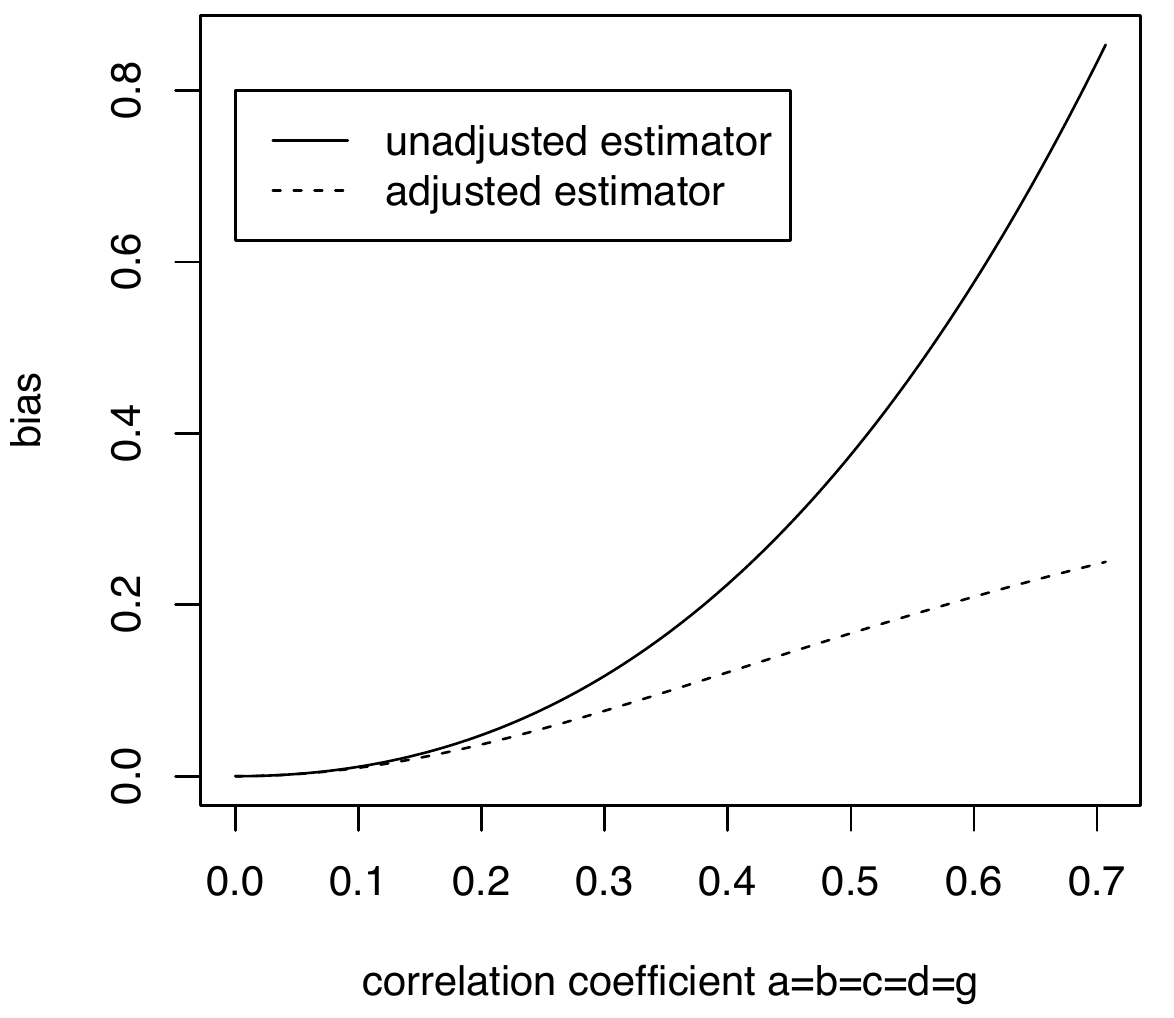} }\\
\\
(a) $M$-Bias, with possible deviations: correlated $U$ and $V$, and an additional arrow from $W$ to $T$
&
(b) Biases of resulting unadjusted and adjusted estimators with $a=b=c=d=g$.
\end{tabular}
\caption{Sensitivity Analysis of \cite{Pearl2009b}'s Critique on \cite{Rubin2007}}\label{fg::rubin-pearl-sensitivity}
\end{figure}

\cite{Pearl2009b} cites \cite{Rubin2007}'s example about the causal effect of smoking habits ($T$) on lung cancer ($Y$), and argues that conditioning on the pretreatment covariate ``seat-belt usage'' ($M$) would introduce spurious associations, since $M$ could be reasonably thought of as an indicator of a person's attitudes toward societal norms ($U$) as well as safety and health related measures ($W$). 
Assuming all the analysis is already conditioned on other observed covariates, we focus our discussion on the five variables $(U, W, M, T, Y)$, of which the dependence structure is illustrated by Figure \ref{fg::rubin-pearl-sensitivity}.
Since the patterns with a continuous treatment and a binary treatment are similar, we focus our discussion on LSEMs.

%
%
%

As \cite{Pearl2009b} points out, 
\begin{quotation}\it
\noindent
If we have good reasons to believe that these two types of attitudes are marginally independent, we have a pure $M$-structure on our hand.
\end{quotation}
In the case with $\rho=0$, conditioning on $M$ will lead to spurious correlation between $T$ and $Y$ under the null, and will bias the estimation of the causal effect of $T$ on $Y$.
However, \cite{Pearl2009b} also recognizes that the independence assumption seems very strong in this example, since $U$ and $W$ are both background variables about the habit and personality of a person.
\cite{Pearl2009b} further argues:
\begin{quotation}\it 
\noindent
But even if marginal independence does not hold precisely, conditioning on ``seat-belt usage'' is likely to introduce spurious associations, hence bias, and should be approached with caution.
\end{quotation}

Although we believe most things should be approached with caution, our work, above, suggests that even mild perturbations of an $M$-Structure can switch which of the two approaches, conditioning or not conditioning, is likely to remove more bias.
In particular,  \cite{Pearl2009b} is correct in that the adjusted estimator indeed tends to introduce more bias than the unadjusted one when an exact $M$-Structure holds and thus the general advice ``to condition on all observed covariates'' may not be sensible in this context. 
However, in the example of \cite{Rubin2007}, the exact independence between a person's attitude toward societal norms $U$ and safety and health related measures $W$ is questionable, since we have good reasons to believe that other hidden variables such as income and family background will affect both $U$ and $W$ simultaneously, and thus Pearl's fears may be unfounded.

To examine this further, we consider two possible deviations from the exact $M$-Structure, and investigate the biases of the unadjusted and adjusted estimators for each.
\begin{enumerate}
[(a)]
\item
(Correlated $U$ and $W$) Assume the DGP follows the DAG in Figure \ref{fg::rubin-pearl-sensitivity}, with an additional correlation between the attitudes $U$ and $W$ as shown in Figure \ref{fg::mbias}.
If we then assume that all the correlation coefficients have the same positive magnitude, earlier results demonstrate that the adjusted estimator is preferable as it strictly dominates the unadjusted estimator except for extremely large values of the correlation coefficients.

Furthermore, in \cite{Rubin2007}'s example, attitudes toward societal norms $U$ are more likely to affect the ``seat-belt usage'' variable $M$ than safety and health related measures $W$, which further strengthens the case for adjustment.
If we were willing to assume that $c$ is zero but $\rho$ is not, equation (\ref{eq::mbias-c0}) in Section \ref{sec::gaussian-m-butterfly} again shows that the adjusted estimator is superior.

\item
(An arrow from $W$ to $T$)
Pearl's example seems a bit confusing on further inspection, even if we accept his independence assumption $U\ind W$.
In particular, one's ``attitudes towards safety and health related measures'' likely impact one's decisions about smoking.
Therefore, we might reasonably expect an arrow from $W$ to $T$.
In Figure \ref{fg::rubin-pearl-sensitivity}(a), we remove the correlation between $U$ and $W$, but we allow an arrow from $W$ to $T$, i.e., the generating equation for $T$ becomes $T = aU + gW + \sqrt{1-a^2-g^2} \varepsilon_T$. Lemma \ref{lemma::rubin-pearl-sensitivity3} in Appendix gives the associated formulae for biases of the adjusted and unadjusted estimators.  Figure \ref{fg::rubin-pearl-sensitivity}(b) shows that, assuming $a=b=c=d=g$ (i.e., equal correlations), the adjusted estimator is uniformly better.

\end{enumerate}

\section{Two Further Issues}\label{sec::further}

One controversy about $M$-Bias is whether $M$-Structure is rare or not in practice, and we go through several examples to discuss this issue. In the second part of this section, we make a distinction between asymptotic and finite sample properties of $M$-Bias.

\paragraph{Is $M$-Structure Rare?}

Although \cite{Pearl2009b} argues that $M$-Bias is a structural property, \cite{Rubin2009} claims that $M$-bias is a rare phenomenon such as ``trying to balance a multidimensional cone on its point with no external supports in some visible directions.''
As mentioned in the introduction, \cite{Gelman2011} argues that, in social sciences, ``true zeros'' are rare and consequently the independence structure in the exact $M$-Structure is also rare. Section \ref{sec::example-rubin-pearl} revisited the controversial example between Professors Pearl and Rubin, and Figure \ref{fg::rubin-pearl-sensitivity}(a) illustrated two possible deviations from the exact $M$-Structure.
Both deviations suggested conditioning is a superior choice.
In the following, we review three other examples of $M$-Bias in the current literature, and investigate the plausibility of the exact $M$-Structure.

\begin{figure}[ht]
\centering
\footnotesize
\begin{tabular}{p{0.33\columnwidth} p{0.33\columnwidth} p{0.33\columnwidth}}
\adjustbox{valign=m}{\includegraphics[width = 0.32\textwidth]{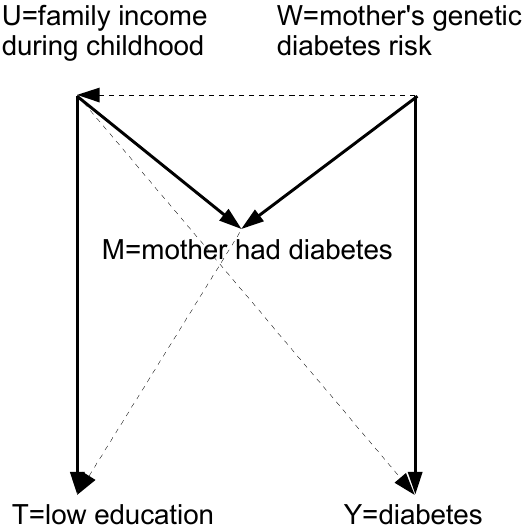}}
&
\adjustbox{valign=m}{\includegraphics[width = 0.32\textwidth]{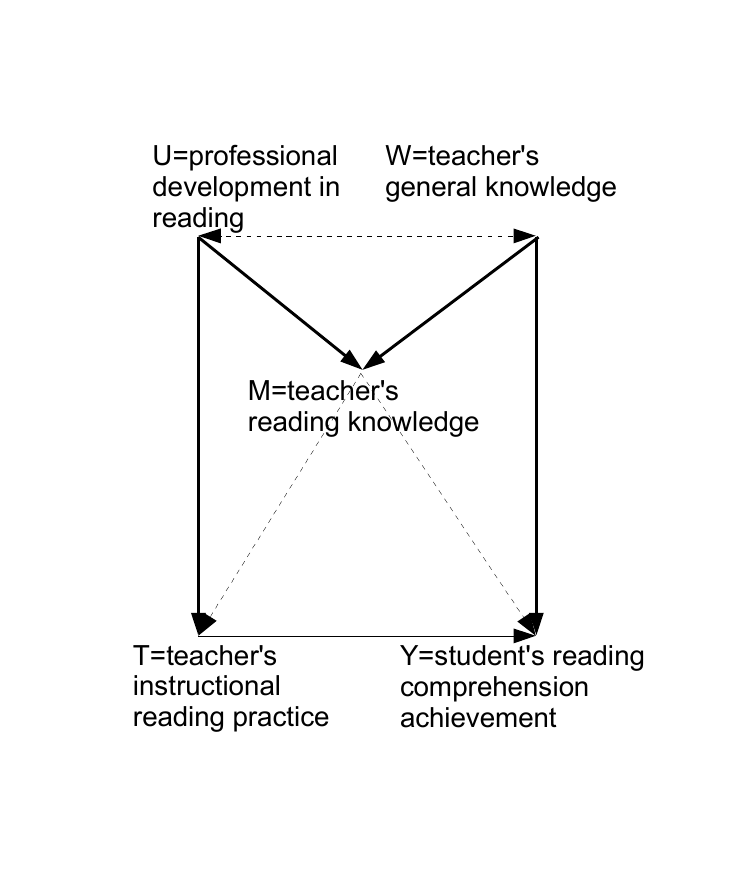} }
&\adjustbox{valign=m}{\includegraphics[width = 0.32\textwidth]{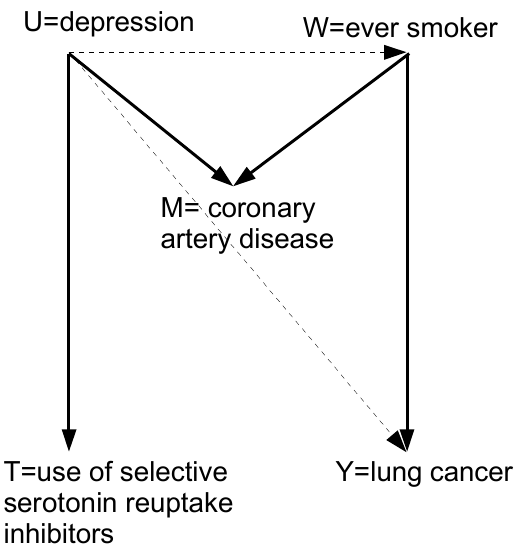} }\\
\\
(a) \cite{Glymour2006}
&
(b) \cite{Kelcey2011} 
&
(c) \cite{Liu2012}
\end{tabular}
\caption{Three $M$-Structures with Possible Deviations in Dashed Arrows}\label{fg::m-rare}
\end{figure}

As shown in Figure \ref{fg::m-rare}(a), \cite{Glymour2006} postulates a possible $M$-Structure with exposure ``low income,'' outcome ``diabetes,'' and $M$ variable ``mother had diabetes,'' where ``family income during childhood'' affects both exposure and $M$, and ``mother's genetic diabetes risk'' affect both outcome and $M$. However, this $M$-Structure is subject to several plausible deviations: ``mother's genetic diabetes risk'' may affect ``family income during childhood; ``mother had diabetes'' may affect ``low education;'' and ``family income during childhood'' may affect ``diabetes.''


\cite{Kelcey2011} have $(T,Y,M,U,W)$ as ``teacher's instructional reading practice,'' ``student's reading comprehension achievement,'' ``teacher's reading knowledge,'' ``professional development in reading,'' and ``teacher's general knowledge.'' See Figure \ref{fg::m-rare}(b). However, this $M$-Structure is dubious because of the possible correlation between the latent $(U, W)$ and the confounding effect of $M$ on the relationship between $T$ and $Y$.

Figure \ref{fg::m-rare}(c) is a possible $M$-Structure investigated by \cite{Liu2012}, where $(T,Y,M,U,W)$ are ``use of selective serotonin reuptake inhibitors (SSRI),'' ``lung cancer,'' ``coronary artery disease,'' ``depression,'' and ``ever smoker.'' Although it is plausible that ``coronary artery disease'' is not the confounder between ``use of SSRI'' and ``lung cancer,'' it is very likely that ``depression'' affects both ``ever smoker'' and ``lung cancer.''

In summary, while all the examples above were quite useful to illustrate $M$-Bias in theoretical research, it is unwise to believe that these $M$-Structures are exact based on our background knowledge. Therefore, we suggest researchers conduct sensitivity analysis, such as illustrated earlier, according to their scientific knowledge about the structure of the DAG and the associated parameters.

\paragraph{Asymptotic versus Finite Sample Properties.}

The discussion in the previous sections are mainly based on asymptotic theory assuming large samples. As argued by \cite{Pearl2009a}, this approach allows for investigating the existence of bias in a certain DAG, and asymptotic analysis helps reveal the structural property of a DAG. A referee pointed out that the asymptotic theory is quite different from the more practical finite sample theory.
In finite sample data analysis, practitioners, often interested in interval estimation and hypothesis testing, are typically more interested in whether associated confidence intervals cover the true causal parameters at nominal rates, and whether tests for null hypotheses about the causal effect have valid size.
These questions are related to the asymptotic property of the DAGs, but also depend on the sample size, the procedure for constructing confidence interval, and choice of test statistic. 
Theoretical discussion of the finite sample theory is unfortunately more difficult. Simulation study, however, is an alternative tool for these questions. Some studies exist. In particular, \cite{Liu2012} simulate large cohort studies under an $M$-Structure corresponding to their science question of interest, and find that the impact of $M$-Bias was small for most of their $178$ scenarios unless the association between $M$ and the unmeasured confounders is very large.

\section{Discussion}\label{sec::discussion}

For objective causal inference, Rubin and Rosenbaum suggest balancing all the pretreatment covariate in observational studies to parallel with the design of randomized experiments \citep{Rubin2007, Rubin2008, Rubin2009, Rosenbaum2002}, which is called the ``pretreatment criterion'' \citep{VanderWeele2011}.
However, Pearl and other researchers \citep{Pearl2009a, Pearl2009b, Shrier2008, Shrier2009, Sjolander2009} criticize the ``pretreatment criterion'' by pointing out that this criterion may lead to biased inference in presence of a possible $M$-Structure even if the treatment assignment is unconfounded.
We investigate this controversy in detail for LSEMs, ideally providing a template for future research about more general DAGs (e.g., nonparametric and nonlinear models). 
While we agree that Pearl's warning is very insightful, our asymptotic theory shows that, at least for LSEMs, this conclusion is quite sensitive to various deviations from the exact $M$-Structure, e.g., to circumstances where latent causes may be correlated or the $M$ variable may also be a confounder between the treatment and the outcome. 
We also go through several candidate $M$-Structures in the existing literature, and find that exact $M$-Structure is likely to be rare with various deviations typically being more plausible. 
Overall, this coupled with our asymptotic theory suggests that for linear systems, except in some extreme cases, adjusting for all the pretreatment covariates is in fact a reasonable choice.

\section*{Acknowledgment}
The authors thank all the participants in the ``Causal Graphs in Low and High Dimensions'' seminar at Harvard Statistics Department in Fall, 2012, and thank Professor Peter Spirtes for sending us his slides \citep{Spirtes2002}. Comments from the associate editor and two reviewers greatly improved the quality of our paper.

\section*{Appendix: Lemmas and Proofs}
\begin{lemma}\label{lemma::2}
In the linear regression model
$
Y = \beta_0 + \beta_T T + \beta_M M + \varepsilon
$
with $\varepsilon\ind (T,M)$ and $E(\varepsilon) = 0$, 
we have
\begin{eqnarray*}
\beta_T &=& \frac{ \Cov(Y, T) \Var(M) - \Cov(Y, M)  \Cov(M, T)   }{ \Var(T) \Var(M) - \Cov^2(M, T)  }.
\end{eqnarray*}
\end{lemma}

\noindent {\it Proof.}
Solve for $(\beta_T, \beta_M)$ using the following moment conditions
$$\left\{
\begin{array}{ccc}
\Cov(Y, T)&=&\beta_T \Var(T) + \beta_M \Cov(M, T),\\
\Cov(Y, M)&=&\beta_T \Cov(M, T) + \beta_M \Var(M).
\end{array}\right.
$$

\begin{lemma}\label{lemma::mbias-gaussian}
Under the model generated by Figure \ref{fg::mbias}, the regression coefficient of $T$ by regressing $Y$ onto $(T,M)$ is
\begin{eqnarray*}
 \beta_T = \frac{  ad\rho(1-b^2-c^2-bc\rho) -  abcd}{  1 - (ab + ac\rho)^2    }.
\end{eqnarray*}
\end{lemma}

\noindent {\it Proof.}
We apply Lemma~\ref{lemma::2}, where all variance terms such as $\Var(M)$ are $1$, and the covariance terms are easily calculated.
For example, we have
$
 \Cov( Y, T ) 
= \Cov( dW, aU ) = ad \rho,
$
and
$
 \Cov( Y, M )  =  \Cov( dW  , bU + cW    ) 
 = bd\rho + cd .
$

\begin{lemma}\label{lemma::butterfly-gaussian}
Under the model generate by Figure \ref{fg::butterfly-bias}, the regression coefficient of $T$ from regressing $Y$ onto $(T,M)$ is
\begin{eqnarray*}
\beta_T =  -  \frac{abcd}{  1 - (ab + e)^2    } . 
\end{eqnarray*}
\end{lemma}

\noindent {\it Proof.} Similar to the proof of Lemma~\ref{lemma::mbias-gaussian}, we apply Lemma \ref{lemma::2}.
%
%

\begin{lemma}\label{lemma::1}
Assume that $(X_1,X_2)$ follows a bivariate Normal distribution with means zero, variances one, and correlation coefficient $r.$
Then 
$
\E(X_1\mid X_2 \geq z) - \E(X_1 \mid X_2< z) = r\eta(z),
$
where $\eta(z) = \phi(z) / \{  \Phi(z) \Phi(-z) \}.$
\end{lemma}

\noindent {\it Proof.}
Since $X_1 = rX_2 + \sqrt{1-r^2} Z$ with $Z\sim N(0,1)$ and $Z\ind X_2$, we have
\begin{eqnarray*}
\E(X_1\mid X_2 \geq z) = r \E(X_2\mid X_2 \geq z) = \frac{r}{\Phi(-z) } \int_z^{\infty} x\phi(x)dx = - \frac{r}{\Phi(-z) } \int_z^{\infty} d \phi(x) = r  \frac{  \phi(z) } { \Phi(-z)} .
\end{eqnarray*}
Similarly, we have $ \E(X_1\mid X_2  <  z) = \E(X_1\mid - X_2  >   -z) = -r\phi(-z)/\Phi(z) = - r\phi(z)/  \Phi( z)$.
Therefore,
$
\E(X_1\mid X_2 \geq z) -\E(X_1\mid X_2  <  z)  = r\phi(z)\left\{  1 / \Phi(-z)   + 1 / \Phi(z)    \right\} = r\eta(z).
$

\begin{lemma}\label{lemma:covariance}
The covariance between $X$ and $B\sim$ Bernoulli$(p)$ is
$$
\Cov(X,B) = p(1-p)\{  \E(X\mid B=1) - \E(X\mid B=0) \}. 
$$
\end{lemma}

\noindent {\it Proof.} It follows from the definition of the covariance.
%

\begin{lemma}\label{lemma::3}
Under the model generated by Figure \ref{fg::binary}(a), the regression coefficient of $T$ from regressing $Y$ onto $(T,M)$ is
\begin{eqnarray*}
\beta_T = 
\frac{ ad\eta(\alpha) \{  \rho (1-b^2-c^2-bc\rho) - bc \}   } { \rho \{     1 - (ab + ac\rho)^2 \phi(\alpha)  \eta(\alpha)   \}   }  . 
\end{eqnarray*}
\end{lemma}

\noindent {\it Proof.}
We have the following joint Normality of $(Y, M, T^*)$:
\begin{eqnarray*}
\begin{pmatrix}
Y\\ M\\ T^*
\end{pmatrix}
&\sim&
\bm{N}_3
\left\{ 
\begin{pmatrix}
0\\ 0\\  0
\end{pmatrix},
\begin{pmatrix}
1     &   bd\rho  + cd  &      ad\rho  \\
bd\rho  + cd & 1 & ab + ac\rho \\
ad\rho &     ab + ac\rho & 1
\end{pmatrix}
\right\}. \label{eq::binaryT-mbias}
\end{eqnarray*}
From Lemma \ref{lemma::1}, we have
\begin{eqnarray*}
\E(M\mid T=1) - \E(M\mid T=0) &=& \E(M  \mid T^* \geq \alpha  ) - \E(M\mid T^*< \alpha ) = (ab + ac\rho) \eta(\alpha),\\
\E(Y \mid T=1) - \E(Y \mid T=0) &=& \E(Y  \mid T^* \geq \alpha  ) - \E(Y \mid T^*< \alpha ) = ad\rho \eta(\alpha).
\end{eqnarray*}
Therefore, from Lemma \ref{lemma:covariance}, the covariances are
$\Cov(M,T) = \Phi(\alpha)\Phi(-\alpha)  (ab + ac\rho) \eta(\alpha),$ and $
\Cov(Y,T)  = \Phi(\alpha)\Phi(-\alpha) ad\rho \eta(\alpha).$
According to Lemma \ref{lemma::2}, the regression coefficient $\beta_T$ is
\begin{eqnarray*}
\beta_T 
&=& \frac{  \Phi(\alpha)\Phi(-\alpha)  ad\rho \eta(\alpha)   - (bd\rho  + cd) \Phi(\alpha)\Phi(-\alpha)  (ab + ac\rho) \eta(\alpha)   }{    \Phi(\alpha)\Phi(-\alpha)  -  \Phi^2(\alpha)\Phi^2(-\alpha)  (ab + ac\rho)^2 \eta^2(\alpha)   }  \\
&=&\frac{ ad\eta(\alpha) \{  \rho (1-b^2-c^2-bc\rho) - bc \}   } { \rho \{     1 - (ab + ac\rho)^2 \phi(\alpha)  \eta(\alpha)   \}   }  .
\end{eqnarray*}

\begin{lemma}\label{lemma::4}
Under the model generated by Figure \ref{fg::binary}(b), the regression coefficient of $T$ from regressing $Y$ onto $(T, M)$ is
\begin{eqnarray*}
\beta_T=   - \frac{  abcd \eta(\alpha)  }{ 1 -  (ab + e)\phi(\alpha)}.
\end{eqnarray*}
\end{lemma}

\noindent {\it Proof.}
We have the following joint Normality of $(Y, M, T^*)$:
\begin{eqnarray*}
\begin{pmatrix}
Y\\ M\\ T^*
\end{pmatrix}
&\sim&
\bm{N}_3
\left\{ 
\begin{pmatrix}
0\\ 0\\  0 
\end{pmatrix},
\begin{pmatrix}
1     &   cd + f  &     cde + abf + ef  \\
cd + f & 1 & ab + e \\
 cde + abf + ef   &     ab + e & 1
\end{pmatrix}
\right\}.\label{eq::binaryT-butterfly}
\end{eqnarray*}

From Lemma \ref{lemma::1}, we have
\begin{eqnarray*}
\E(M\mid T=1) - \E(M\mid T=0) &=& \E(M  \mid T^* \geq \alpha  ) - \E(M\mid T^*< \alpha ) = (ab + e) \eta(\alpha),\\
\E(Y \mid T=1) - \E(Y \mid T=0) &=& \E(Y  \mid T^* \geq \alpha  ) - \E(Y \mid T^*< \alpha ) = (cde+abf+ef) \eta(\alpha).
\end{eqnarray*}
From Lemma \ref{lemma:covariance}, we obtain their covariances
$\Cov(M,T)  =  \Phi(\alpha) \Phi(-\alpha) (ab + e) \eta(\alpha),$ and 
$\Cov(Y,T)   =  \Phi(\alpha) \Phi(-\alpha) (cde + abf + ef)\eta(\alpha).$
According to Lemma \ref{lemma::2}, the regression coefficient $\beta_T$ is
\begin{eqnarray*}
\beta_T 
&=& \frac{  \Phi(\alpha)\Phi(-\alpha)  (cde + abf + ef)  \eta(\alpha)   - (cd + f ) \Phi(\alpha)\Phi(-\alpha)  (ab + e) \eta(\alpha)   }{    \Phi(\alpha)\Phi(-\alpha)  -  \Phi^2(\alpha)\Phi^2(-\alpha)  (ab + e)^2 \eta^2(\alpha)   } \\
&=&  - \frac{  abcd \eta(\alpha)  }{ 1 -  (ab + e)^2 \phi(\alpha) \eta(\alpha)}.
\end{eqnarray*}

\begin{lemma}
\label{lemma::rubin-pearl-sensitivity3}
Under the model generated by Figure \ref{fg::rubin-pearl-sensitivity}(a) with an arrow from $W$ to $T$, the unadjusted estimator has bias $ad\rho + dg$, and the adjusted estimator has bias
$$
\frac{   dg  - (cd) (  ab + cg )   }{1 - ( ab + cg   )^2 }.
$$
\end{lemma}

\noindent {\it Proof.}
The unadjusted estimator is  
$
\Cov(T, Y) = ad\rho + dg.
$
Expanding Lemma \ref{lemma::2} gives the above as the regression coefficient of $T$ from regressing $Y$ onto $(T, M)$.

\bibliographystyle{DeGruyter}
\bibliography{JCI_Mbias}

\end{document}